\documentclass{article}

\usepackage{silence}
\usepackage[english]{babel}
\usepackage[utf8]{inputenc}
\usepackage[T1]{fontenc}

\usepackage[scale=0.8]{geometry}

\usepackage{graphicx}
\usepackage{xcolor}

\usepackage{booktabs}
\usepackage{enumitem}
\usepackage{makecell}

\usepackage{amssymb}
\usepackage{mathtools}
\usepackage{faktor}
\usepackage{cases}

\usepackage[amsmath,amsthm,hyperref,thmmarks]{ntheorem}

\usepackage{pgfplots} 
\usepackage{pgfplotstable}
\usepackage{fp}
\pgfplotsset{compat=newest}
\usetikzlibrary{
	fpu,
	fixedpointarithmetic,
	babel,
	external,
	arrows,
	plotmarks,
	positioning,
	angles,
	quotes,
	intersections,
	calc,
	spy,
	decorations.pathreplacing,
	matrix,
	fit,
}
\usepgfplotslibrary{fillbetween}
\usepackage{shellesc}
\tikzset{
	external/system call={
			lualatex \tikzexternalcheckshellescape -halt-on-error -interaction=batchmode -jobname "\image" "\texsource"
		}
}
\makeatletter  \def\pgfplotscolormappdfmax@inv{1000} \makeatother

\makeatletter

\newcommand{\printslope}[4]{
	\def\nero@printslope@orderlist{#1}
	\edef\nero@printslope@xpos{#2}
	\edef\nero@printslope@ypos{#3}
	\edef\nero@printslope@width{#4}
	\pgfmathparse{\nero@printslope@xpos+\nero@printslope@width}
	\edef\nero@printslope@px{\pgfmathresult}
	\edef\nero@printslope@py{\nero@printslope@ypos}
	\edef\nero@printslope@qx{\nero@printslope@xpos}
	\edef\nero@printslope@ry{\nero@printslope@ypos}
	\foreach \nero@printslope@order in {#1}{
			\pgfmathparse{
				((\nero@printslope@px/\nero@printslope@xpos)^(\nero@printslope@order/2))*\nero@printslope@ypos}
			\edef\nero@printslope@qy{\pgfmathresult}
			\edef\nero@aux1{\noexpand\draw[line width=0.7pt]
				(axis cs:\nero@printslope@xpos,\nero@printslope@ry)
				-- (axis cs:\nero@printslope@qx,\nero@printslope@qy)
				-- (axis cs:\nero@printslope@px,\nero@printslope@py);}
			\nero@aux1
			\pgfmathparse{10^((ln(\nero@printslope@ry)+ln(\nero@printslope@qy))/(ln(10)*2))}
			\edef\nero@printslope@labelpos{\pgfmathresult}
			\edef\nero@aux2{\noexpand\node[anchor=east] at
				(axis cs:\nero@printslope@qx,\nero@printslope@labelpos)
				{\noexpand\small \nero@printslope@order};}
			\nero@aux2
			\global\edef\nero@printslope@ry{\nero@printslope@qy}
		}
	\draw[line width=0.7pt] (axis cs:\nero@printslope@xpos,\nero@printslope@ypos)
	|- (axis cs:\nero@printslope@px,\nero@printslope@py);
	\pgfmathparse{10^((ln(\nero@printslope@px)+ln(\nero@printslope@xpos))/(ln(10)*2))}
	\edef\nero@printslope@labelpos{\pgfmathresult}
	\node[anchor=north] at (axis cs:\nero@printslope@labelpos,\nero@printslope@ypos) {\small 1};
}

\makeatother

\pgfplotsset{
	colormap={rainbow desaturated}{
			rgb=(0.278431, 0.278431, 0.858824)
			rgb=(0, 0, 0.360784  )
			rgb=(0, 1, 1  )
			rgb=(0, 0.501961, 0  )
			rgb=(1, 1 , 0  )
			rgb=(1, 0.380392 , 0  )
			rgb=(0.419608, 0  , 0 )
			rgb=(0.878431, 0.301961  , 0.301961 )
		}
}

\pgfplotsset{
	/pgfplots/surf shading/precision=pdf,
}

\pgfplotstableset{
	font={\small},
	columns/N/.style  = {sci zerofill,   precision=2, column name = \makecell{$N$}},
	columns/h_HSR_o1_error/.style = {sci zerofill,   precision=2, column name = \makecell{error}},
	columns/h_HSR_o1_order/.style = {fixed zerofill, precision=2, column name = \makecell{order}},
	columns/h_HDR_o1_error/.style = {sci zerofill,   precision=2, column name = \makecell{error}},
	columns/h_HDR_o1_order/.style = {fixed zerofill, precision=2, column name = \makecell{order}},
	columns/h_HDR_o2theta_error/.style = {sci zerofill,   precision=2, column name = \makecell{error}},
	columns/h_HDR_o2theta_order/.style = {fixed zerofill, precision=2, column name = \makecell{order}},
	columns/h_HDR_o3theta_error/.style = {sci zerofill,   precision=2, column name = \makecell{error}},
	columns/h_HDR_o3theta_order/.style = {fixed zerofill, precision=2, column name = \makecell{order}},
	columns/HSR_o1/.style = {sci, sci zerofill, precision=2, column name = \makecell{HSR, $\mathbb{P}_0$}},
	columns/HDR_o1/.style = {sci, sci zerofill, precision=2, column name = \makecell{HDR, $\mathbb{P}_0$}},
	columns/HDR_o2theta/.style = {sci, sci zerofill, precision=2, column name = \makecell{HDR, $\mathbb{P}_1$}},
	columns/HDR_o3theta/.style = {sci, sci zerofill, precision=2, column name = \makecell{HDR, $\mathbb{P}_2$}},
	columns/variable/.style = {string type, column name = \space, postproc cell content/.code={%
					\pgfplotsutilstrreplace{\_}{ }{##1}%
					\pgfkeyslet{/pgfplots/table/@cell content}\pgfplotsretval
				},postproc cell content/.style={@cell content=error on $##1$},},
	col sep=space,
	empty cells with={---},
	every head row/.style={before row=\toprule,after row=\midrule},
	every last row/.style={after row=\bottomrule},
	set thousands separator={},
}

\definecolor{fill_topo}{RGB}{191,191,191}

\newcommand{\lp}{\left(}
\newcommand{\rp}{\right)}
\newcommand{\lb}{\left[}
\newcommand{\rb}{\right]}
\newcommand{\R}{\mathbb{R}}
\newcommand{\Z}{\mathbb{Z}}

\newcommand{\bp}{\begin{pmatrix}}
		\newcommand{\ep}{\end{pmatrix}}
\newcommand{\dx}{{\Delta x}}

\newcommand{\dt}{{\Delta t}}

\newcommand{\sgn}{\mathop{\mathrm{sgn}}}                      
\newcommand{\cst}{\mathop{\mathrm{cst}}}                      
\newcommand{\sfr}[2]{{\scriptstyle \faktor {#1} {#2}}}        

\newcommand{\px}{\partial_x}

\newcommand{\cO}{\mathcal{O}}

\newcommand{\dZ}{{\Delta Z}}
\newcommand{\im}{{i-1}}

\newcommand{\imhm}{{i-\frac{1}{2}, -}}
\newcommand{\imh}{{i-\frac{1}{2}}}
\newcommand{\imhp}{{i-\frac{1}{2}, +}}
\newcommand{\ipmhm}{{i\pm\frac{1}{2}, -}}
\newcommand{\ipmh}{{i\pm\frac{1}{2}}}
\newcommand{\ipmhp}{{i\pm\frac{1}{2}, +}}
\newcommand{\iphm}{{i+\frac{1}{2}, -}}
\newcommand{\iph}{{i+\frac{1}{2}}}
\newcommand{\iphp}{{i+\frac{1}{2}, +}}
\newcommand{\iphpm}{{i+\frac{1}{2}, \pm}}

\newcommand{\ip}{{i+1}}

\DeclareMathOperator{\Fr}{{Fr}}
\DeclareMathOperator{\cF}{{\mathcal{F}}}
\DeclareMathOperator{\cH}{{\mathcal{H}}}
\DeclareMathOperator{\cS}{{\mathcal{S}}}
\providecommand{\abs}[1]{\lvert#1\rvert}

\numberwithin{equation}{section}

\makeatletter
\newsavebox\myboxA
\newsavebox\myboxB
\newlength\mylenA

\newcommand*\xoverline[2][0.75]{%
	\sbox{\myboxA}{$\m@th#2$}%
	\setbox\myboxB\null
	\ht\myboxB=\ht\myboxA%
	\dp\myboxB=\dp\myboxA%
	\wd\myboxB=#1\wd\myboxA
	\sbox\myboxB{$\m@th\overline{\copy\myboxB}$}
	\setlength\mylenA{\the\wd\myboxA}
	\addtolength\mylenA{-\the\wd\myboxB}%
	\ifdim\wd\myboxB<\wd\myboxA%
		\rlap{\hskip 0.5\mylenA\usebox\myboxB}{\usebox\myboxA}%
	\else
		\hskip -0.5\mylenA\rlap{\usebox\myboxA}{\hskip 0.5\mylenA\usebox\myboxB}%
	\fi}

\makeatother

\makeatletter
\renewcommand*\env@matrix[1][\arraystretch]{%
	\edef\arraystretch{#1}%
	\hskip -\arraycolsep
	\let\@ifnextchar\new@ifnextchar
	\array{*\c@MaxMatrixCols c}}
\makeatother

\usepackage[colorlinks=true]{hyperref}
\usepackage[noabbrev,capitalise,nameinlink]{cleveref}

\hypersetup{colorlinks=true,
	pdfauthor = {Christophe Berthon, Victor Michel-Dansac},
	linkcolor=blue,
	citecolor=blue,
	filecolor=blue,
	urlcolor=blue
}

\newtheorem{theorem}{Theorem}
\newtheorem{lemma}[theorem]{Lemma}
\newtheorem{proposition}[theorem]{Proposition}

\theoremstyle{remark}

\theoremstyle{definition}

\renewcommand{\dh}{{\Delta h}}

\makeatletter
\newcommand{\vast}{\bBigg@{4}}
\newcommand{\Vast}{\bBigg@{5}}
\makeatother

\title{A fully well-balanced hydrodynamic reconstruction}

\author{
    Christophe Berthon\thanks{
        Université de Nantes, CNRS UMR 6629, Laboratoire de Mathématiques Jean Leray, 2 rue de la Houssinière, BP 92208, 44322 Nantes, France;
        email address: christophe.berthon@univ-nantes.fr;
        homepage: \url{https://www.math.sciences.univ-nantes.fr/~berthon}
    } \and
    Victor Michel-Dansac\thanks{
        Université de Strasbourg, CNRS, Inria, IRMA, F-67000 Strasbourg, France;
        email address: victor.michel-dansac@inria.fr;
        homepage: \url{https://irma.math.unistra.fr/~micheldansac}
    }~\thanks{corresponding author}
}

\date{April 29, 2024}

\usepackage{xcolor}
\newcommand{\ra}{\color{black}}

\begin{document}

\maketitle

\paragraph*{Abstract}
The present work focuses on the numerical approximation of the weak solutions
of the shallow water model over a non-flat topography.
In particular, we pay close attention to steady solutions with nonzero velocity.
The goal of this work is to derive a scheme that
exactly preserves these stationary solutions,
as well as the commonly preserved lake at rest steady solution.
    {\ra These moving steady states are solution to a nonlinear equation.
        We emphasize that the method proposed here never requires solving this nonlinear equation;
        instead, a suitable linearization is derived.}
To address this issue,
we propose an extension of the well-known hydrostatic reconstruction.
By appropriately defining the reconstructed states at the interfaces,
any numerical flux function, combined with a relevant source term discretization,
produces a well-balanced scheme that preserves both moving and non-moving steady solutions.
This eliminates the need to construct specific numerical fluxes.
Additionally, we prove that the resulting scheme is consistent with
the homogeneous system on flat topographies,
and that it reduces to the hydrostatic reconstruction when the velocity vanishes.
To increase the accuracy of the simulations,
{\ra we propose a well-balanced high-order procedure,
        which still does not require solving any nonlinear equation}.
Several numerical experiments demonstrate the effectiveness of the numerical scheme.

\section{Introduction}

In this work, we consider the numerical approximation of the weak
solutions of the shallow water equations, given as follows:
\begin{equation}
	\label{eq:shallow_water}
	\begin{dcases}
		\partial_t h + \partial_x q = 0, \\
		\partial_t q + \partial_x \lp \frac {q^2} h + \frac 1 2 g h^2 \rp = - g h \partial_x Z,
	\end{dcases}
\end{equation}
where $h \geq 0$ is the water height and $q \in \mathbb{R}$ stands for
the discharge.  The given topography function $Z$ is assumed to be smooth
enough, and $g > 0$ is the gravity constant. For the sake of
convenience in the notation, as long as~$h>0$, we define the velocity
$u$ as follows:
\begin{equation}
	\label{eq:def_velocity}
	u = \frac q h.
\end{equation}
In fact, defining the velocity in dry regions, when the water height $h$ vanishes,
requires special attention.
In this paper, we impose the following conventions:
\begin{equation}\label{eq:limh0}
	\lim_{h\to 0} \frac{q}{h}=0
	\text{\quad and \quad}
	\frac{q^2}{h^3}
	\text{ is bounded when } h \to 0.
\end{equation}
The first convention is commonly used when simulating dry areas.
The second convention, however, is less usual.
It is associated with the behavior of the Froude number in dry areas.

For the sake of clarity, we introduce the following condensed notation
\begin{equation}
	\label{eq:def_W}
	W = \bp h \\ q \ep
	\text{, \quad}
	F(W) = \bp q \vphantom{\dfrac 1 2} \\ \dfrac {q^2} h + \dfrac 1 2 g h^2 \ep
	\text{\quad and \quad}
	S(W) = \bp 0 \\ - g h \partial_x Z \ep,
\end{equation}
so that the shallow water model \eqref{eq:shallow_water} rewrites
\begin{equation}
	\label{eq:shallow_water_condensed}
	\partial_t W + \partial_x F(W) = S(W).
\end{equation}
This system is endowed with an initial data
\begin{equation}
	\label{eq:initial_data}
	W(x, t=0) = W^0(x),
\end{equation}
where $W^0$ is defined according to prescribed physics.

When deriving numerical schemes to approximate the weak solutions of \eqref{eq:shallow_water},
special attention must be given to the steady solutions, see \cite{BerVaz1994,GreLeR1996}.
Such solutions of major interest are governed by the system $\partial_x F(W) = S(W)$.
Solving this time-independent system yields smooth steady solutions that only depend on $x$,
given by
\begin{equation}
	\label{eq:steady_with_constants}
	\begin{dcases}
		q \eqqcolon q_0, \\
		\frac {q^2} {2 h^2} + g (h + Z) \eqqcolon B_0,
	\end{dcases}
\end{equation}
where $q_0$ and $B_0$ are given constant values and with a Froude
number $\Fr^2=\frac{q^2}{gh^3}$ different from $1$.

Over the last three decades,
a significant amount of research has been devoted to the derivation of numerical schemes
able to exactly preserve the steady solutions.
After pioneering work by Berm\'udez and V\'azquez \cite{BerVaz1994},
and next by Greenberg and Leroux \cite{GreLeR1996},
it is now well-known that steady solutions play a crucial role
when designing numerical discretizations.
Indeed, even when simulating time-dependent phenomena,
small errors in the approximation of a steady solution
may accumulate over time and render the simulation irrelevant.
To avoid such an issue, well-balanced schemes have been
introduced~(see~\cite{GreLeR1996}, as well as \cite{BerVaz1994} for the C-property).
A scheme is said to be well-balanced if
it exactly preserves at least some steady solutions,
and we emphasize that, in general, a naive discretization of \eqref{eq:shallow_water}
will not lead to a well-balanced scheme.
In \cite{Gos2000}, Gosse obtains a well-balanced scheme by
solving the nonlinear Bernoulli equation
\eqref{eq:steady_with_constants} in each cell and at each time step.
Next, in \cite{Jin2001}, Jin derives an approximately well-balanced scheme,
where accurately approximating the steady solutions,
rather than exactly preserving them,
results in good scheme behavior.

However, these approaches have several drawbacks.
Indeed, exactly solving Bernoulli's equation is computationally costly,
and the approximately well-balanced property may be insufficient in some cases.
An interesting alternative was introduced by Audusse et al.
in \cite{AudBouBriKlePer2004},
where the \emph{hydrostatic reconstruction} was designed to
exactly preserve the lake at rest steady state.
The lake at rest is governed by particular case of~\eqref{eq:steady_with_constants},
where~$q_0$ is fixed equal to $0$, thus leading to
$h + Z \eqqcolon H_0$, with $H_0$ a given constant free surface.
Base on the simplicity of this steaady solution, Audusse et al. derived a technique
able to ensure that a given numerical scheme exactly preserves the lake at rest.
The simplicity of this approach makes it very attractive.

Indeed, let us briefly recall the hydrostatic reconstruction.
We introduce a space discretization made of cells~$(x_{i - 1/2}, x_{i + 1/2})$
of constant size $\Delta x>0$, such that $x_{i + 1/2}=x_{i - 1/2}+\Delta x$. We
also introduce the time discretization $t^{n+1}=t^n+\Delta t$, where
the time increment $\Delta t>0$ is restricted by a CFL-type
condition. To define the hydrostatic reconstruction, the authors
of \cite{AudBouBriKlePer2004} merely have to adopt a numerical flux
function $\cF \lp W_i^n, W_{i+1}^n \rp$ that is consistent with the exact flux
function $F(W)$ given by \eqref{eq:def_W}, i.e.
$\cF(W,W)=F(W)$. Next, at each interface~$x_\iph$, they introduce
the following reconstructed water heights:
\begin{equation*}
    h_{\iph,-}^n = \max \lp 0,h_i^n+Z_i-Z_\iph \rp
    \text{\quad and \quad}
    h_{\iph,+}^n = \max \lp 0,h_{i+1}^n+Z_{i+1}-Z_\iph \rp,
\end{equation*}
with $Z_\iph=\max(Z_i,Z_{i+1})$, where
\begin{equation*}
	Z_i=\frac{1}{\Delta x} \int_{x_\imh}^{x_\iph} Z(x) dx
\end{equation*}
is a finite volume discretization of the topography function $Z$.
Equipped with such a reconstruction,
with $W_i^n$ a given approximation of $W(x,t^n)$ over the
cell $(x_\imh,x_\iph)$, the approximate solution of
\eqref{eq:shallow_water} at time $t^{n+1}$ is given by
\begin{equation}
	\label{eq:scheme}
	W_i^{n+1} =
	W_i^n
	- \frac \dt \dx \lp \cF \lp W_\iphm^n, W_\iphp^n \rp - \cF \lp W_\imhm^n, W_\imhp^n \rp \rp
	+ \dt \cS_i^n,
\end{equation}
where
\begin{equation*}
	W_\iphm^n = \begin{pmatrix}
		h_\iphm^n \\ h_\iphm^n u_i^n
	\end{pmatrix},
	\quad
	W_\iphp^n = \begin{pmatrix}
		h_\iphp^n \\ h_\iphp^n u_{i+1}^n
	\end{pmatrix}
	\text{\quad and \quad}
	u_i^n = \frac{(hu)_i^n}{h_i^n}.
\end{equation*}
Concerning the source term $\cS_i^n$, it now defined as follows:
\begin{equation*}
	\cS_i^n = \begin{pmatrix}
		0 \\
		\dfrac g 2 \dfrac{(h_\iphm^n)^2 - (h_\imhp^n)^2}{\Delta x}
	\end{pmatrix}.
\end{equation*}

Thanks to its simplicity and versatility, the hydrostatic reconstruction is frequently used
(for instance, see~\cite{MorCasPar2013,NatRibTwa2014,CheNoe2017,BerDurFouSalZab2019},
for applications and extensions)
to get a scheme that exactly preserves, specifically, the lake at rest steady solution.
Several other techniques have also been developed to yield lake at rest-preserving schemes
(for instance, see the non-exhaustive list
\cite{XinZhaShu2010,XinShu2011,BerFou2012,Delestre2013,CheNoe2017,Li2018,DurVilBar2020}).

Yet, merely considering the lake at rest preservation may be insufficient in some simulations,
see \cite{XinShuNoe2011}.
As a consequence, more recently,
new approaches have been proposed to deal with still or moving steady states,
which are described by \eqref{eq:steady_with_constants} with $q = 0$ or $q \neq 0$.
For a few examples of schemes that exactly preserve moving steady states,
the reader is referred
to \cite{NoeXinShu2007,Xin2014} for high-order schemes,
to \cite{BouMor2010} for a Suliciu relaxation scheme dealing with the subcritical case,
to \cite{BerCha2016,MicBerClaFou2016,MicBerClaFou2017,MicBerClaFou2021}
for Godunov-type schemes that exactly capture moving steady states,
or to~\cite{GomCasPar2021} for a control-based approach.

In this work, we are interested in an extension of
the hydrostatic reconstruction \cite{AudBouBriKlePer2004}
in order to exactly capture moving steady states.
	{\ra Contrary to other techniques (e.g. \cite{Gos2000,CasPar2020}),
		we emphasize that the proposed reconstruction does not
		require solving the nonlinear equations \eqref{eq:steady_with_constants}.}
This issue is addressed by designing new reconstructed states
$W_{i+1/2,\pm}^n$ at the interface $x_{i+1/2}$.
Unlike the generalized hydrostatic reconstruction from \cite{CasParPar2007},
where the nonlinear Bernoulli equation in solved in each cell
and where the positivity of the water height may be lost,
we here introduce a suitable linearization.
The relevant properties of the hydrostatic reconstruction,
namely the ability to deal with transitions between wet and dry areas,
as well as its versatility,
are also satisfied by the technique derived in this paper.

To address such issues, the paper is structured as follows.
Firstly, \cref{sec:comments_steady} contains a few necessary comments
related to the steady states \eqref{eq:steady_with_constants}. Indeed,
the smoothness of the steady
solutions \eqref{eq:steady_with_constants} is a priori governed by the
smoothness of the topography function $Z$. However, because of the
topography discretization, the required smoothness is lost. As a
consequence, a discontinuous extension
of \eqref{eq:steady_with_constants} must be considered. In fact, the
discontinuous extension of the steady solution turns out to be the
main ingredient of the forthcoming hydrodynamic extension.
Then, in \cref{sec:scheme}, we introduce an easy interface reconstruction technique
that preserves the moving and non-moving steady states,
building on the versatility of
the hydrostatic reconstruction~\cite{AudBouBriKlePer2004}.
Moreover, we prove that the derived hydrodynamic reconstruction
can handle transitions between wet and dry areas.
Since the hydrodynamic reconstruction is defined according to a
sequence of properties, in \cref{sec:choice_of_H}, we specify
the reconstruction, and we present an explicit definition.
Next, in \cref{sec:HO_WB}, we suggest {\ra an easy} high-order accurate
technique which preserves the moving and non-moving steady states.
Finally, in \cref{sec:numerics}, we present numerical experiments
that illustrate the relevance of the hydrodynamic reconstruction.

\section{Discretization of steady solutions}
\label{sec:comments_steady}

This section is devoted to some comments about the discretization of
the smooth steady solutions given by \eqref{eq:steady_with_constants}.
First, let us emphasize that arguing the smoothness of the topography function $Z(x)$,
we have $Z_\ip - Z_i = \cO(\dx)$.

Now, we consider a non-dry region; namely with $h_i^n>0$. For a given
pair $(q_0,B_0)\in\R^2$ and an arbitrary time~$t^n$,
according to \eqref{eq:steady_with_constants} a discrete steady
solution is naturally given by, for all $i\in\Z$,
\begin{equation}
	\label{eq:steady_discrete}
	\begin{dcases}
		q_i^n \eqqcolon q_0, \\
		\frac {(q_i^n)^2} {2 (h_i^n)^2} + g (h_i^n + Z_i) \eqqcolon B_0.
	\end{dcases}
\end{equation}
In fact, we now show that this natural choice of the approximate
steady solutions may introduce inconsistencies coming from a loss of
smoothness in \eqref{eq:steady_with_constants}. Indeed,
from \eqref{eq:steady_discrete}, we easily obtain a local per interface
definition of a steady solution, given as follows for all $i\in\Z$:
\begin{equation}
	\label{eq:steadyperinterface}
	\begin{dcases}
		q_i^n  = q_\ip^n, \\
		\frac {(q_i^n)^2} {2 (h_i^n)^2} + g (h_i^n + Z_i) =
		\frac {(q_\ip^n)^2} {2 (h_\ip^n)^2} + g (h_\ip^n + Z_\ip).
	\end{dcases}
\end{equation}
Since, at the interface located at $x_{i+1/2}$, the pair $(Z_i,Z_\ip)$
defines a discontinuous topography function with a small
jump in $\cO(\dx)$, the local definition \eqref{eq:steadyperinterface}
is not sufficient to ensure the smoothness required to derive~\eqref{eq:steady_with_constants}.
Actually, even if the error to the smoothness is small and controlled
by $\cO(\dx)$, this failure may create non-physical steady solutions.

Indeed, let us momentarily consider the equivalent augmented
system \cite{LeFTha2007} as follows:
\begin{equation}
	\label{eq:shallow_water_augmented}
	\begin{dcases}
		\partial_t h + \partial_x q = 0,    \\
		\partial_t q + \partial_x \lp \frac {q^2} h + \frac 1
		2 g h^2 \rp + g h \partial_x Z = 0, \\
		\partial_t Z = 0.
	\end{dcases}
\end{equation}
We easily see that this system contains a stationary wave endowed with
the following Riemann invariants:
\begin{equation}
	\label{eq:steadyIR}
	q
	\text{\quad and \quad}
	\frac{q^2}{2h^2}+g(h+Z).
\end{equation}
Naturally, from these Riemann invariants, we easily recover the steady
solutions \eqref{eq:steady_with_constants}. As a
consequence, \eqref{eq:steadyperinterface} exactly coincides with the
preservation of the Riemann invariants for the stationary wave.
In other words,~\eqref{eq:steadyperinterface} is relevant in defining
discrete steady solutions independently of the smoothness of the
pair $(Z_i,Z_\ip)$.

However, it is essential to note that \eqref{eq:steadyperinterface}
can produce physically inconsistent solutions. Indeed,
as demonstrated in~\cite{LeFTha2007}, the Riemann problem of
the shallow water equations \eqref{eq:shallow_water}
-- or equivalently \eqref{eq:shallow_water_augmented} --
with discontinuous topography may admit non-unique solutions.
While the steady solution preserving the Riemann invariants~\eqref{eq:steadyIR}
is a possible solution,
it is not the only one and does not seem to be
the expected physical solution \cite{LeFTha2007,Aleksyuk2022}.
{\ra Put in other words, let us consider a discontinuous initial condition,
with discontinuous topography,
made of two states that satisfy constant Riemann invariants \eqref{eq:steadyIR}.
After~\cite{LeFTha2007,LeFTha2011}, such a Riemann problem has multiple solutions,
one made of a stationary contact wave,
and others corresponding to solutions
containing non-stationary shock and rarefaction waves.
At this level, we conjecture that the stationary contact wave solution
is numerically unstable,
while at least the another one seems numerically stable.}
Such an assertion will be
illustrated by numerical experiments displayed in \cref{sec:numerics}.

Another important comment about the steady solutions concerns their
definition in the case of a partially dry space domain.
Indeed, after~\cite{MicBerClaFou2016},
{\ra moving steady solutions and dry regions cannot coexist}.
As a consequence, as soon as a dry area is present, smooth steady
solutions must be at rest.

\section{A family of moving steady states-preserving schemes}
\label{sec:scheme}

In this section, we present a simple state interface reconstruction method
that enables the numerical method to preserve both moving and non-moving steady solutions
as given by \eqref{eq:steadyperinterface},
following the approach by Audusse et al.~\cite{AudBouBriKlePer2004}.
With a scheme given by \eqref{eq:scheme},
in the spirit of the usual hydrostatic reconstruction,
we here design the reconstructed states
$W_{{i+1/2},\pm}^n$ at the interface $x_{i+1/2}$ and the source term
discretization $\cS_i^n$ in order to preserve the expected moving
steady solutions \eqref{eq:steady_discrete}.

The interface reconstructed states are now given by
$W_\iphpm^n = (h_\iphpm^n, q_\iphpm^n)^\intercal$,
where we have set
\begin{equation}
	\label{eq:hydrodynamic_reconstruction}
	\begin{dcases}
		h^n_\iphm = \max\left(0,h^n_i
		+ \lp Z_i - Z_\iph \rp
		+ 2 \Fr^2(h^n_i, h^n_\iph, q^n_i) \cH \lp h^n_i, h^n_\iph, q^n_i, Z_\iph - Z_i \rp\right),           \\
		h_\iphp^n = \max\left(0,h^n_\ip
		+ \lp Z_\ip - Z_\iph \rp
		+ 2 \Fr^2(h^n_\ip, h^n_\iph, q^n_\ip) \cH \lp h^n_\ip, h^n_\iph, q^n_\ip, Z_\iph - Z_\ip \rp\right), \\
		q_\iphm^n = q_i^n,                                                                                   \\
		q_\iphp = q_\ip^n,
	\end{dcases}
\end{equation}
with an intermediate reconstruction $(W_\iph, Z_\iph)$ at the
interface $x_\iph$ defined by
\begin{equation}
	\label{eq:intermediate_reconstruction}
	(W_\iph^n, Z_\iph) =
	\begin{dcases}
		(W_i^n, Z_i)     & \text{if } Z_i > Z_\ip, \\
		(W_\ip^n, Z_\ip) & \text{otherwise.}
	\end{dcases}
\end{equation}
Moreover, we have introduced the approximate Froude number as follows:
\begin{equation}
	\label{eq:Fr2}
	\Fr^2(h^n_i, h^n_\iph, q^n_i) = \frac{
		(q_i^n)^2(h_i^n+h_\iph^n)}{2g(h_i^n)^2(h_\iph^n)^2}.
\end{equation}
Next, concerning the source term discretization,
$\cS_i^n = (0, (\cS_q)_i^n)^\intercal$,
we have adopted the following definition:
\begin{equation}
	\begin{aligned}
		\dx (\cS_q)_i^n =
		 & - g \frac {2 h_\imhp^n h_\iphm^n} {h_\imhp^n + h_\iphm^n} \lp Z_\iph - Z_\imh \rp                  \\
		 & + \frac {4 g} {h_\imhp^n + h_\iphm^n} \cH \lp h_\imhp^n, h_\iphm^n, q_i^n, Z_\iph - Z_\imh \rp^3.
	\end{aligned}
	\label{eq:source}
\end{equation}

At this level, it is worth noticing that the hydrodynamic reconstruction
cannot be fully characterized until the function $\cH$ is defined.
A possible choice of the particular function is
detailed in \cref{sec:choice_of_H}.

From now on, let us emphasize that the interface
reconstruction \eqref{eq:hydrodynamic_reconstruction} is nothing but
the standard hydrostatic reconstruction
from \cite{AudBouBriKlePer2004}, augmented with an additional term
governed by the function $\cH$.  Of course, this new term serves to perturb the
hydrostatic reconstruction in order to provide a hydrodynamic
reconstruction that is capable of preserving the moving steady state solutions.

Now, we impose suitable hypotheses to be satisfied by the perturbation
$\cH$ so that the hydrodynamic reconstruction scheme \eqref{eq:scheme} is consistent,
preserves the steady solutions \eqref{eq:steady_discrete},
and efficiently handles wet/dry transitions.

To address this issue, the function $\cH$ must be endowed with several properties.
\begin{enumerate}[label=($\cH$-\arabic*)]
	\item \label{item:def:H_property_dZ_0}
	      In order to recover the required consistency, the perturbation
	      $\cH:\R^\star_+ \times \R^\star_+ \times \R \times \R \to \R$
	      must be continuous. Moreover, since the moving steady solutions
	      are non-constant only for non-flat topographies,
	      the perturbation $\cH$ should vanish in the case of a flat topography.
	      Hence, $\cH$ must satisfy the following asymptotic behavior
	      for all $h_L>0$, $h_R>0$ and $\bar q\in\R$:
	      \begin{equation*}
		      \lim_{\Delta Z\to 0} \cH(h_L,h_R,\bar q,\Delta Z) =0.
	      \end{equation*}
	      Note that, here, $h_L$ represents either $h_i^n$ or $h_\ip^n$,
	      $h_R$ represents $h_\iph^n$, and
	      $\bar q$ represents either $q_i^n$ or~$q_\ip^n$.

	\item \label{item:def:H_property_WB}
	      Next, the well-balanced property relies on a technical condition,
	      whose relevance will become apparent in a forthcoming proof.
	      This property will be recovered by imposing
	      \begin{equation*}
		      \cH \lp h_L, h_R, \bar q, \dZ \rp = \frac{1}{2}(h_R-h_L),
	      \end{equation*}
	      for all $(h_L,h_R,\bar q,\Delta
		      Z)\in \R^\star_+ \times \R^\star_+ \times \R \times \R$ such
	      that $\dZ = - \lp h_R - h_L \rp \lp 1 - \Fr^2(h_L, h_R, \bar q) \rp$,
	      and with $\Fr^2(h_L, h_R, q_0) \neq 1$ so that $h_R-h_L=\cO(\Delta Z)$.

	\item \label{item:def:properties_of_H_dry_wet}
	      The last property we enforce concerns the wet/dry transition.
	      To that end, we have to make sure that all the involved
	      quantities are well-defined in dry regions.
	      To address such an issue, we require that the following limit holds
	      in order for $h_\imhp^n$ and~$h_\iphm^n$ to be bounded:
	      \begin{equation*}
		      \lim_{h_i^n\to0}
		      \Fr^2(h^n_i, h^n_{i\pm\frac{1}{2}}, q^n_i)
		      \cH \lp h_i^n, h_{i\pm\frac{1}{2}}^n, q_i^n, Z_{i\pm\frac{1}{2}} -
		      Z_i \rp
		      =0.
	      \end{equation*}
	      Next, to make sure that the source term vanishes in dry areas,
	      we also have to impose
	      \begin{equation*}
		      \lim_{h_i^n\to0}
		      \frac{1}{h_\imhp^n + h_\iphm^n} \cH \lp h_\imhp^n, h_\iphm^n, q_i^n, Z_\iph - Z_\imh \rp^3 =0
		      \text{\qquad and \qquad}
		      \lim_{h_i^n\to0}
		      \frac {h_\imhp^n h_\iphm^n} {h_\imhp^n + h_\iphm^n} =0.
	      \end{equation*}

\end{enumerate}

Before we state the main properties satisfied by the scheme
\eqref{eq:scheme}--\eqref{eq:hydrodynamic_reconstruction}--\eqref{eq:source},
let us underline once again that $\cH$ is nothing but a small perturbation.
This is due to the fact that the topography function is assumed to be smooth,
resulting in $Z_{i+1}-Z_i$ being of order $\Delta x$.
As a consequence of \ref{item:def:H_property_dZ_0}, we have
	{\ra $\cH(h_i^n,h_{i+1/2}^n,q_i^n,Z_{i+1/2}-Z_i)=o(\Delta x)$
		and
		$\cH(h_\ip^n,h_{i+1/2}^n,q_\ip^n,Z_{i+1/2}-Z_\ip)=o(\Delta x)$}
for all
$i\in\Z$. In this sense, the application $\cH$ is clearly a small
perturbation of the original hydrostatic reconstruction designed
in \cite{AudBouBriKlePer2004}. Moreover, it is also worth mentioning that obtaining an
explicit definition of $\cH$, such that
assumptions \ref{item:def:H_property_dZ_0}, \ref{item:def:H_property_WB}
and \ref{item:def:properties_of_H_dry_wet} hold, is not a trivial task.
\cref{sec:choice_of_H} is devoted to exhibiting an admissible
perturbation. More specifically, we will show that the set of admissible
perturbation according to assumptions \ref{item:def:H_property_dZ_0}, \ref{item:def:H_property_WB}
and \ref{item:def:properties_of_H_dry_wet} is not empty, and we will design
suitable approximations of the admissible functions $\cH$.

In addition, let us emphasize that
assumptions \ref{item:def:properties_of_H_dry_wet} are necessary to
define a wet/dry transition.
In dry regions, where $h_{i-1}^n=h_i^n=h_{i+1}^n=0$,
the hydrodynamic reconstruction \eqref{eq:hydrodynamic_reconstruction}
and the source term \eqref{eq:source} are well-defined,
and they vanish due to conventions \eqref{eq:limh0}.

We are now able to state our main result,
which outlines the properties of the scheme \eqref{eq:scheme}
endowed with the hydrodynamic reconstruction \eqref{eq:hydrodynamic_reconstruction}
and the source term discretization \eqref{eq:source}.

\begin{theorem}
	\label{thm:scheme_properties}
	Let $\cH(h_L, h_R, q_0, \dZ)$ be a function which satisfies the
	assumptions  \ref{item:def:H_property_dZ_0}, \ref{item:def:H_property_WB}
	and \ref{item:def:properties_of_H_dry_wet}.
	For non-negative water heights $h_i^n\ge0$ for all $i\in\Z$,
	the scheme \eqref{eq:scheme}--\eqref{eq:hydrodynamic_reconstruction}--\eqref{eq:source} satisfies the following properties:
	\begin{enumerate}[label=(\thetheorem-\alph*)]
		\item \label{item:thm:scheme_properties_consistency}
		      it is consistent with the shallow water equations \eqref{eq:shallow_water};

		\item \label{item:thm:scheme_properties_WB}
		      it preserves the steady states with nonzero velocity,
		      in the sense that if $(W_i^n)_{i\in\Z}$
		      satisfy \eqref{eq:steadyperinterface} for all $i\in\Z$ and
		      $\Fr^2(h_i^n,h_{i\pm\frac{1}{2}}^n,q_i^n)\not=1$ then
		      $W_i^{n+1}=W_i^n$;

		\item \label{item:thm:scheme_properties_dry_wet_positivity}
		      it is non-negativity-preserving,
		      i.e., if $h_i^n \geq 0$, then $h_i^{n+1} \geq 0$.

	\end{enumerate}
\end{theorem}

Now, in order to establish this main result, we first need to prove the following two lemmas.
They state some needed properties to be satisfied by the hydrodynamic reconstruction.
The first lemma deals with wet areas, while the second one addresses the wet/dry and dry cases.

\begin{lemma}
	\label{lem:hydrodynamic_reconstruction}
	For positive water heights $h_i^n>0$ for all $i\in\Z$ far away from dry areas,
	the hydrodynamic reconstruction \eqref{eq:hydrodynamic_reconstruction},
	with assumptions \ref{item:def:H_property_dZ_0} and \ref{item:def:H_property_WB},
	satisfies the following properties:
	\begin{enumerate}[label=(\thetheorem-\alph*)]
		\item \label{item:lem:hydrodynamic_reconstruction_flat_topo}
		      if $Z_i = Z_\ip$, then $h_\iphm^n = h_i^n$ and $h_\iphp^n = h_\ip^n$;
		\item \label{item:lem:hydrodynamic_reconstruction_static}
		      the hydrodynamic reconstruction
		      \eqref{eq:hydrodynamic_reconstruction} degenerates to the standard hydrostatic reconstruction from \cite{AudBouBriKlePer2004}
		      as soon as $q_i^n = q_\ip^n = 0$;
		\item \label{item:lem:hydrodynamic_reconstruction_WB}
		      if two consecutive states $(W_i^n, Z_i)$ and
		      $(W_\ip^n, Z_\ip)$ satisfy the local per interface steady state
		      definition~\eqref{eq:steadyperinterface},
		      then $h_\iphm^n = h_\iph^n$ and $h_\iphp^n = h_\iph^n$.
	\end{enumerate}
\end{lemma}

\begin{proof}
	The proof of \ref{item:lem:hydrodynamic_reconstruction_flat_topo} is immediate,
	since it relies on property \ref{item:def:H_property_dZ_0}.

	Similarly, by inspection, we note that if $q_i^n = q_\ip^n = 0$
	in \eqref{eq:hydrodynamic_reconstruction},
	then $\Fr^2(h_i, h_\ipmh, 0) = 0$, to get
	\begin{equation*}
		\label{eq:hydrodynamic_reconstruction_for_vanishing_q}
		\begin{dcases}
			h_\iphm^n = \max\left(0,h_i^n
			+ \lp Z_i - Z_\iph \rp\right), \\
			h_\iphp^n = \max\left(0,h_\ip^n
			+ \lp Z_\ip - Z_\iph \rp\right),
		\end{dcases}
	\end{equation*}
	which is nothing but the standard hydrostatic reconstruction,
	and thus \ref{item:lem:hydrodynamic_reconstruction_static} holds.

	Regarding \ref{item:lem:hydrodynamic_reconstruction_WB}, since
	$(W_i, Z_i)$ and $(W_\ip, Z_\ip)$ define a local per interface
	steady state according to \eqref{eq:steadyperinterface},
	then we have
	\begin{equation*}
		\label{eq:definition_B0}
		\begin{dcases}
			q_i^n = q_\ip^n \eqqcolon q_0, \\
			\frac {q_0^2} {2 (h_i^n)^2} + g (h_i^n + Z_i) = \frac {q_0^2} {2 (h_\ip^n)^2} + g (h_\ip^n + Z_\ip) \eqqcolon B_0.
		\end{dcases}
	\end{equation*}
	Next, by definition of $(h_\iph, Z_\iph)$ given by \eqref{eq:intermediate_reconstruction},
	we get
	\begin{equation*}
		\label{eq:B0_applied_to_interface}
		\frac {q_0^2} {2 (h_\iph^n)^2} + g (h_\iph^n + Z_\iph) = B_0.
	\end{equation*}
    As a consequence, we have
    \begin{equation*}
        \frac{q_0^2}{2(h_i^n)^2} + g(h_i^n+Z_i) =
        \frac{q_0^2}{2(h_\iph^n)^2} + g(h_\iph^n+Z_\iph)
        \text{\quad and \quad}
        \frac{q_0^2}{2(h_\ip^n)^2} + g(h_\ip^n+Z_\ip) =
        \frac{q_0^2}{2(h_\iph^n)^2} + g(h_\iph^n+Z_\iph),
    \end{equation*}
    to get
    \begin{equation*}
        Z_\iph-Z_i = - (h_\iph^n-h_i^n)\Big( 1-\Fr^2(h_i^n,h_\iph^n,q_i^n)\Big)
        \text{\quad and \quad}
        Z_\iph-Z_\ip = - (h_\iph^n-h_\ip^n)\Big( 1-\Fr^2(h_\ip^n,h_\iph^n,q_\ip^n)\Big).
    \end{equation*}
    Next, arguing property \ref{item:def:H_property_WB}, we obtain
    \begin{equation*}
        \cH \lp h_i^n, h_\iph^n, q_i^n, Z_\iph - Z_i \rp = \frac 1 2 \lp h_\iph^n - h_i^n \rp
        \text{\quad and \quad}
        \cH \lp h_\ip^n, h_\iph^n, q_\ip^n, Z_\iph - Z_\ip \rp = \frac 1 2 \lp h_\iph^n - h_\ip^n \rp.
    \end{equation*}
	Plugging these values into the definition \eqref{eq:hydrodynamic_reconstruction}
	yields the following chain of equalities:
	\begin{equation*}
		\label{eq:h_iphm_is_h_iph}
		\begin{aligned}
			h_\iphm ^n
			 & =
			\max\left(0,
			h_i^n
			+ \lp Z_i - Z_\iph \rp
			+ \frac{(q_i^n)^2 (h_i^n + h_\iph^n)}{g (h_i^n)^2 (h_\iph^n)^2} \frac 1 2 \lp h_\iph^n - h_i^n \rp
			\right),                          \\
			 & =
			\max\left(0,
			h_i^n
			+ \lp Z_i - Z_\iph \rp
			+ \frac{q_0^2}{2 g} \lp \frac 1 {(h_i^n)^2}
			- \frac 1 {(h_\iph^n)^2} \rp
			\right),                          \\
			 & =
			\max\left(0,
			\frac {B_0} g - Z_\iph - \frac{q_0^2}{2 g} \frac 1 {(h_\iph^n)^2}
			\right),                          \\
			 & = \max\left(0,h_\iph^n\right), \\
			 & = h_\iph^n.
		\end{aligned}
	\end{equation*}
	Similar relations lead to $h_\iphp = h_\iph$,
	which concludes the proof of \cref{lem:hydrodynamic_reconstruction}.
\end{proof}

Next, in addition to the above result, it is also necessary to establish the
behavior of the hydrodynamic reconstruction~\eqref{eq:hydrodynamic_reconstruction} in dry areas.
This is the object of the following lemma.

\begin{lemma}
	\label{prop:properties_of_HSR_dry_wet}
	Assume $\cH$ satisfies the assumptions \ref{item:def:H_property_dZ_0}, \ref{item:def:H_property_WB}
	and \ref{item:def:properties_of_H_dry_wet}.
	Then, the hydrodynamic
	reconstruction~\eqref{eq:hydrodynamic_reconstruction} verifies:
	\begin{enumerate}[label=(\thetheorem-\alph*)]
		\item \label{item:DW_req_1}
		      if $h_i^n = 0$, then $h_\imhp^n = h_\iphm^n = 0$;
		\item \label{item:DW_req_2}
		      if $h_i^n = 0$, then $\cS_i^n = 0$.
	\end{enumerate}
\end{lemma}

\begin{proof}
	Let us first focus on property \ref{item:DW_req_1}. Because of
	assumptions \ref{item:def:properties_of_H_dry_wet}, arguing
	the water height reconstruction~\eqref{eq:hydrodynamic_reconstruction}
	in the cell	$(x_{i-1/2},x_{i+1/2})$ with $h_i^n=0$, we immediately get
	\begin{equation*}
		h_\imhp^n=\max(0,Z_{i-1}-Z_\imh)
		\text{\quad and \quad}
		h_\iphm^n=\max(0,Z_i-Z_\iph).
	\end{equation*}
	Since we have $Z_\iph=\max(Z_i,Z_{i+1})$ for all $i\in\Z$,
	we then obtain $h_\imhp^n = 0$ and $h_\iphm^n = 0$.

	Next, concerning \ref{item:DW_req_2}, once again
	assumptions \ref{item:def:properties_of_H_dry_wet} immediately
	enforce $\cS_i^n = 0$ as soon as $h_i^n=0$, and the proof is achieved.
\end{proof}

From now on, it should be noted that properties \ref{item:DW_req_1} and \ref{item:DW_req_2},
which comprise \cref{prop:properties_of_HSR_dry_wet},
are satisfied by the standard hydrostatic reconstruction from \cite{AudBouBriKlePer2004}.

Using the intermediate results established in \cref{lem:hydrodynamic_reconstruction} and
\cref{prop:properties_of_HSR_dry_wet}, we can now proceed to prove \cref{thm:scheme_properties}.

\begin{proof}[Proof of \cref{thm:scheme_properties}]

	Regarding the consistency in \ref{item:thm:scheme_properties_consistency},
	let us recall that the numerical flux $\cF$ is assumed to be
	consistent with the exact flux function $F(W)$ given by \eqref{eq:def_W}.
	Since this consistent flux function $\cF$
	is applied without modification to the reconstructed states
	given by \eqref{eq:hydrodynamic_reconstruction},
	the consistency of the flux is maintained,
	based on the proof given in \cite{AudBouBriKlePer2004}.
	We still need to prove the consistency of the source term \eqref{eq:source}.
	More specifically, we have to show that $(\cS_q)_i^n$ is consistent with $- g h \px Z$.
	The approximate source term reads
	\begin{equation}
		\label{eq:source_over_dx_consistency}
		(\cS_q)_i^n
		=
		- g \frac {2 h_\imhp^n h_\iphm^n} {h_\imhp^n + h_\iphm^n} \frac {Z_\iph - Z_\imh} \dx
		+ \frac {4 g} {h_\imhp^n + h_\iphm^n} \frac {\cH \lp h_\imhp^n, h_\iphm^n, q_i^n, Z_\iph - Z_\imh \rp^3} \dx.
	\end{equation}
	Arguing \ref{item:def:H_property_dZ_0}, we obtain that the right part of $(\cS_q)_i^n$
	vanishes when $\dx$ approaches $0$ due to the smoothness assumption on the topography.
	Then, it is obvious that the remainder is consistent with $- g h \px Z$,
	which concludes the proof of~\ref{item:thm:scheme_properties_consistency}.

	Next, let us establish the well-balanced property \ref{item:thm:scheme_properties_WB}.
	Assume that $(W^n_\im, Z_\im)$, $(W^n_i, Z_i)$ and $(W^n_\ip, Z_\ip)$
	define the same steady state \eqref{eq:steady_discrete}, with
	constant discharge $q_0$ and Bernoulli's constant $B_0$.
	We now have to prove that $(W^{n+1}_i, Z_i) = (W^n_i, Z_i)$.

	According to \cref{lem:hydrodynamic_reconstruction} (property \ref{item:lem:hydrodynamic_reconstruction_WB}),
	the hydrodynamic reconstruction \eqref{eq:hydrodynamic_reconstruction} becomes
	\begin{equation*}
		\label{eq:hydrodynamic_reconstruction_when_WB}
		\begin{dcases}
			h_\ipmhm^n = h_\ipmhp^n = h_\ipmh^n, \\
			q_\ipmhm^n = q_\ipmhp^n = q_\ipmh^n = q_0.
		\end{dcases}
	\end{equation*}
	Let us set
	$W_{i\pm\frac{1}{2}}^n=(h_{i\pm\frac{1}{2}}^n,q_0)^\intercal$ so that
	the scheme \eqref{eq:scheme} now reads
	\begin{equation}
		\label{eq:scheme_when_WB}
		W_i^{n+1} =
		W_i^n
		- \frac \dt \dx \lp \cF \lp W_\iph^n, W_\iph^n \rp - \cF \lp W_\imh^n, W_\imh^n \rp \rp
		+ \dt \cS_i^n.
	\end{equation}
	Since $\cF$ is a consistent numerical flux, we know that $\cF(W, W) = F(W)$,
	where $F$ is the flux of the shallow water equations defined
	by \eqref{eq:def_W}. Thus, \eqref{eq:scheme_when_WB} writes
	\begin{equation}
		\label{eq:scheme_when_WB_vector}
		\begin{pmatrix}
			h_i^{n+1} \\
			q_i^{n+1}
		\end{pmatrix}
		=
		\begin{pmatrix}
			h_i^n \\
			q_i^n
		\end{pmatrix}
		- \frac \dt \dx
		\begin{pmatrix}
			q_\iph^n - q_\imh^n \\
			\dfrac {(q_\iph^n)^2} {h_\iph^n} + \dfrac 1 2 g (h_\iph^n)^2
			- \dfrac {(q_\imh^n)^2} {h_\imh^n} - \dfrac 1 2 g (h_\imh^n)^2
		\end{pmatrix}
		+ \dt
		\begin{pmatrix}
			0 \\
			(\cS_q)_i^n
		\end{pmatrix}.
	\end{equation}
	We immediately note that $W_i^{n+1} = W_i^n=(h_i^n,q_0)^\intercal$ for all $i \in \mathbb{Z}$
	as soon as the source term $(\cS_q)_i^n$ satisfies
	\begin{equation}
		\label{eq:source_when_steady}
		\dx (\cS_q)_i^n
		=
		\dfrac{q_0^2}{h_\iph^n} + \dfrac 1 2 g (h_\iph^n)^2
		-
		\dfrac{q_0^2}{h_\imh^n} - \dfrac 1 2 g (h_\imh^n)^2.
	\end{equation}
	Now, to establish the above relation, since we have
	$h_{i+1/2,\pm}^n=h_{i+1/2}^n$ and $h_{i-1/2,\pm}^n=h_{i-1/2}^n$, the
	source term~$(\cS_q)_i^n$ defined
	by \eqref{eq:source_over_dx_consistency} reads
	\begin{equation}
		\label{eq:Sqaux}
			(\cS_q)_i^n
			=
			 - g \frac {2 h_\imh^n h_\iph^n} {h_\imh^n + h_\iph^n} \frac {Z_\iph - Z_\imh} \dx
			 + \frac {4 g} {h_\imh^n + h_\iph^n} \frac {\cH \lp h_\imh^n, h_\iph^n, q_0, Z_\iph - Z_\imh \rp^3} \dx.
	\end{equation}
	Next, we remark that the states $(W_{i-1/2}^n,Z_{i-1/2})$ and
	$(W_{i+1/2}^n,Z_{i+1/2})$ satisfy the local per interface steady state
	condition \eqref{eq:steadyperinterface}. We then have
	\begin{equation}
		\label{eq:deltaZaux}
		Z_\iph-Z_\imh =
		-(h_\iph^n - h_\imh^n)
		\Big( 1- \Fr^2(h_\imh^n,h_\iph^n,q_0)\Big),
	\end{equation}
	and, as a consequence of hypothesis \ref{item:def:H_property_WB},
	we get
	\begin{equation*}
		\cH(h_\imh^n,h_\iph^n,q_0,Z_\iph-Z_\imh) = \frac{1}{2}(h_\iph^n - h_\imh^n).
	\end{equation*}
	Therefore, \eqref{eq:Sqaux} now reads
	\begin{equation}
		\label{eq:source_successive_expressions_steady_1}
		\dx (\cS_q)_i^n
		=
		- g \frac {2 h_\imh^n h_\iph^n} {h_\imh^n + h_\iph^n} \lp Z_\iph - Z_\imh \rp
		+ g \frac {\lp h_\iph^n - h_\imh^n \rp^3} {2 \lp h_\imh^n + h_\iph^n \rp}.
	\end{equation}
	Plugging \eqref{eq:deltaZaux} into \eqref{eq:source_successive_expressions_steady_1},
	the above relation reformulates
	\begin{equation}
		\label{eq:source_successive_expressions_steady_2}
		\dx (\cS_q)_i^n  =
		\dfrac{q_0^2}{h_\iph^n} + \dfrac 1 2 g (h_\iph^n)^2
		-
		\dfrac{q_0^2}{h_\imh^n} - \dfrac 1 2 g (h_\imh^n)^2,
	\end{equation}
	which concludes the proof of the well-balanced
	property \ref{item:thm:scheme_properties_WB}.

	The final property \ref{item:thm:scheme_properties_dry_wet_positivity}
	is a direct consequence of the definition of $h_{i+1/2,\pm}^n$,
	which involves taking a maximum with $0$.
	Then, to prove the non-negativity preservation
	satisfied by the scheme, we exactly follow the proof given in \cite{AudBouBriKlePer2004}.
	The establishment of \cref{thm:scheme_properties} is thus completed.
\end{proof}

\section{\texorpdfstring{One possible choice for the function $\cH$}{One possible choice for the function H}}
\label{sec:choice_of_H}

The goal of this section is to propose an expression for the function $\cH$ that
satisfies the required properties~\ref{item:def:H_property_dZ_0}
through \ref{item:def:properties_of_H_dry_wet}.
Recall that $\cH$ is a function from~$\R^\star_+ \times \R^\star_+ \times \R \times \R$ to $\R$,
applied for instance on~$(h^n_i, h^n_{i+1/2}, q^n_i, Z_{i+1/2} - Z_i)$
to compute $h_{i \pm 1/2, -}$.
For the sake of clarity, throughout this section,
$\cH$ will be written as a function $\cH(h_L, h_R, \bar q, \dZ)$.
Its arguments may be omitted for more concise notation.

As a first step, consider $\cH$ given as a solution to the following polynomial equation of degree five, in the spirit of~\cite{BerCha2016}:
\begin{equation}
	\label{eq:H_as_polynomial_solution}
	2 \cH \lp g \lp \bar h^2 - \cH^2 \rp^2 - \bar q^2 \bar h \rp = - g \dZ \lp \bar h^2 - \cH^2 \rp^2,
\end{equation}
where $\bar h = (h_L + h_R) / 2$.
This expression allows us to state the following result.

\begin{lemma}
    \label{lem:H_as_polynomial_solution}
	If the correct solution $\cH$ of the polynomial equation
    \eqref{eq:H_as_polynomial_solution} is chosen,
	then $\cH$ satisfies the consistency property \ref{item:def:H_property_dZ_0}
	and the well-balanced property \ref{item:def:H_property_WB}.
\end{lemma}

\begin{proof}
	For this proof, we consider data far from a dry area.

	Let us start with property \ref{item:def:H_property_WB}.
	To that end, we assume a steady solution with $\Fr^2(h_L, h_R, \bar q) \neq 1$, and we take
	$\dZ = - \lp h_R - h_L \rp \lp 1 - \Fr^2(h_L, h_R, \bar q) \rp$.
	Given the expression \eqref{eq:Fr2} of the squared Froude number,
	plugging the above value of $\dZ$ in \eqref{eq:H_as_polynomial_solution} leads to:
	\begin{equation}
		\label{eq:H_as_polynomial_solution_with_SS}
		2 \cH \lp g \lp \bar h^2 - \cH^2 \rp^2 - \bar q^2 \bar h \rp
		-
		(h_R - h_L) \lp g \lp \bar h^2 - \cH^2 \rp^2 - \bar q^2 \bar h \frac{\lp \bar h^2 - \cH^2 \rp^2}{h_L^2 h_R^2} \rp = 0.
	\end{equation}
	Now, we prove that $\cH = (h_R - h_L) / 2$
	is a solution to the fifth-degree equation
	\eqref{eq:H_as_polynomial_solution_with_SS}.
	First, let us note that
	\begin{equation*}
		\label{eq:h_bar_minus_H_is_hl_hr}
		\frac 1 {h_L h_R} \lp \bar h^2 - \frac{(h_R - h_L)^2} {4} \rp
		=
		\frac 1 {h_L h_R} \frac{(h_R + h_L)^2 - (h_R - h_L)^2} {4}
		=
		1.
	\end{equation*}
	Then, plugging $\cH = (h_R - h_L) / 2$ in the left-hand side of
	\eqref{eq:H_as_polynomial_solution_with_SS} yields:
	\begin{equation*}
		\label{eq:nonlinear_H_is_hr_minus_hl}
            2 \frac{(h_R - h_L)} {2}
            \lp g \lp \bar h^2 - \frac{(h_R - h_L)^2} {4} \rp^2 - \bar q^2 \bar h \rp
            -
            (h_R - h_L)
            \lp g \lp \bar h^2 - \frac{(h_R - h_L)^2} {4} \rp^2 - \bar q^2 \bar h \rp
            =
            0,
	\end{equation*}
	which proves that $\cH = (h_R - h_L) / 2$
	is one of the multiple solutions to \eqref{eq:H_as_polynomial_solution_with_SS}.
    This proves property \ref{item:def:H_property_WB}.

	We turn to property \ref{item:def:H_property_dZ_0}.
    The solution to \eqref{eq:H_as_polynomial_solution_with_SS} we have just exhibited
    satisfies property \ref{item:def:H_property_WB},
    i.e. it satisfies $\cH = (h_R - h_L) / 2 = \mathcal{O}(h_R - h_L)$
    as soon as $\dZ = - \lp h_R - h_L \rp \lp 1 - \Fr^2(h_L, h_R, \bar q) \rp$ for
    $\Fr^2(h_L, h_R, \bar q) \neq 1$.
    However, note that, when $\Delta Z = 0$, this condition is also satisfied,
    and we get, in this case, $h_R - h_L = \mathcal{O}(\dZ)$.
    This immediately proves that the previously exhibited solution $\cH$ satisfies
    $\cH \smash{\underset{\dZ \to 0}{=}} \cO(\dZ)$,
    thereby proving property \ref{item:def:H_property_dZ_0}.

    The proof is thus concluded, since we have exhibited a solution
    to \eqref{eq:H_as_polynomial_solution_with_SS} that satisfies
    both properties \ref{item:def:H_property_dZ_0} and \ref{item:def:H_property_WB}.
\end{proof}

Therefore, according to \cref{lem:H_as_polynomial_solution},
the nonlinear equation \eqref{eq:H_as_polynomial_solution}
has a solution that satisfies both
properties~\ref{item:def:H_property_dZ_0} and \ref{item:def:H_property_WB}.
However,
finding the correct solution is a complex process
that would negate all the benefits of our {\ra linearized} approach.
Indeed, we would need to use Newton's method at each interface
and each time step to compute the solutions to \eqref{eq:H_as_polynomial_solution}.
Moreover, we would then have to choose the correct solution
from among (at most) five possible ones.
Thus, we elect not to pursue this nonlinear direction.
Instead, we provide a relevant linearization of part of
\eqref{eq:H_as_polynomial_solution}.

Assuming that $\cH \neq \pm \bar h$, \eqref{eq:H_as_polynomial_solution} rewrites
\begin{equation}
	\label{eq:H_as_polynomial_solution_with_division}
	2 \cH \lp 1 - \frac{\bar q^2 \bar h}{\lp \bar h^2 - \cH^2 \rp^2} \rp = - \dZ.
\end{equation}
We temporarily assume that $h_L \neq h_R$,
and we set $\dh = h_R - h_L$.
We suggest the following linearization around~$\cH = \Delta h / 2$
of the expression in brackets in \eqref{eq:H_as_polynomial_solution_with_division},
thus modifying the equation satisfied by $\cH$,
to get a quadratic equation in $\cH$:
\begin{equation*}
	\label{eq:H_as_quadratic_solution_without_Froude}
    2 \cH \lp
        1 - \frac{\bar q^2 (h_L + h_R)}{2 g h_L^2 h_R^2}
        + 4 \sgn(\dZ) \sqrt{\frac{\abs{\dZ}}{\abs{\dh}^3}} \lp \dh - 2 \cH \rp
    \rp
    =
    - \dZ,
\end{equation*}
This expression can be simplified by remarking that
\begin{equation*}
    \label{eq:Froude_to_simplify_quadratic}
    1 - \Fr^2(h_L, h_R, \bar q) = 1 - \frac{\bar q^2 (h_L + h_R)}{2 g h_L^2 h_R^2},
\end{equation*}
to get
\begin{equation*}
	\label{eq:H_as_quadratic_solution}
    2 \cH \lp
        1 - \Fr^2 + 4 \sgn(\dZ) \sqrt{\frac{\abs{\dZ}}{\abs{\dh}^3}} \lp \dh - 2 \cH \rp
    \rp
    =
    - \dZ.
\end{equation*}
Solving this quadratic equation for $\cH$ leads to
\begin{equation}
	\label{eq:solution_of_quadratic_H}
    \cH = \frac 1 4 \vast( \dh + \frac {1 - \Fr^2} 4 \sgn(\dZ) \sqrt{\frac{\abs{\dh}^3}{\abs{\dZ}}}
    \pm \sqrt{\lp \dh + \frac {1 - \Fr^2} 4 \sgn(\dZ) \sqrt{\frac{\abs{\dh}^3}{\abs{\dZ}}} \rp^2 + \sqrt{\abs{\dZ} \abs{\dh}^3}} \vast).
\end{equation}
Choosing the correct sign for the $\pm$ in \eqref{eq:solution_of_quadratic_H}
makes it possible to state, and prove, the following result.

\begin{proposition}
    \label{prop:expression_of_H}
	Let
	\begin{equation}
		\label{eq:cH_version_2}
		\cH  = \frac 1 4 \lp E - \sgn(1 - \Fr^2) \sgn(\dZ)
		\sqrt{E^2 + \sqrt{\abs{\dZ} \abs{\dh}^3}} \rp,
        \text{\quad with \quad}
        E = \dh + \frac {1 - \Fr^2} 4 \sgn(\dZ) \sqrt{\frac{\abs{\dh}^3}{\abs{\dZ}}}.
	\end{equation}
	Then $\cH$ satisfies the required properties
	\ref{item:def:H_property_dZ_0}
	(if $\dZ$ and $1 - \Fr^2$ do not simultaneously vanish),
	\ref{item:def:H_property_WB} and
	\ref{item:def:properties_of_H_dry_wet}.
\end{proposition}

\begin{proof}

We prove the three properties
\ref{item:def:H_property_dZ_0},
\ref{item:def:H_property_WB} and
\ref{item:def:properties_of_H_dry_wet} in order.

\begin{enumerate}[label=($\cH$-\arabic*)]

	\item Let us first note that \eqref{eq:cH_version_2} contains a division by $\dZ$.
	      Nevertheless, this expression turns out to be infinitely continuously differentiable around $\dZ = 0$.

	      To prove \ref{item:def:H_property_dZ_0},
	      let us compute the limit of $\cH$ when $\dZ$ goes to $0$.
	      To that end, for the sake of simplicity,
	      we consider the case where $\dZ > 0$, $\dh > 0$ and $1 - \Fr^2 > 0$.
	      In this case, the Taylor expansion provided in
          \cref{sec:Taylor_expansions_H}
          proves \ref{item:def:H_property_dZ_0}.
	      An immediate consequence is that $\cH$ tends to $0$ as $\dZ$ tends to $0$,
	      despite the \emph{a priori} indeterminate division by $\sqrt{\abs{\dZ}}$.
	      Investigating the other cases ($1 - \Fr^2 < 0$ or $\dh \leq 0$)
	      yields the same limit.

	\item We now prove that the expression of $\cH$ satisfies property \ref{item:def:H_property_WB},
	      i.e. that $\cH = \dh / 2$ when a steady solution has been reached.
	      We therefore assume that the solution is steady, i.e. that $\dZ = - (1 - \Fr^2) \dh$,
	      so that the expression of $E$ is simplified as follows:
	      \begin{equation*}
		      \label{eq:E_v2_steady}
              \begin{aligned}
		      E
		      & = \dh + \frac {1 - \Fr^2} 4 \sgn(-(1 - \Fr^2) \dh)
		      \sqrt{\frac{\abs{\dh}^3}{\abs{(1 - \Fr^2)} \abs{\dh}}} \\
		      & = \dh \lp 1 - \frac 1 4 \sqrt{\abs{1 - \Fr^2}} \rp,
              \end{aligned}
	      \end{equation*}
	      and $E^2 + \sqrt{\abs{\dZ} \abs{\dh}^3}$ becomes:
	      \begin{equation*}
		      \label{eq:sqrt_v2_steady}
		      E^2 + \sqrt{\abs{\dZ} \abs{\dh}^3}
		      = \abs{\dh}^2 \lp 1 + \frac 1 4 \sqrt{\abs{1 - \Fr^2}} \rp^2.
	      \end{equation*}
	      Plugging the two expressions above into $\cH$, we get the following chain of equalities:
	      \begin{equation*}
		      \label{eq:H_v2_steady}
		      \begin{aligned}
			      \cH
			       & = \frac 1 4 \lp
			      \dh \lp 1 - \frac 1 4 \sqrt{\abs{1 - \Fr^2}} \rp
			      + \sgn(\dh)
			      \sqrt{\abs{\dh}^2 \lp 1 + \frac 1 4 \sqrt{\abs{1 - \Fr^2}} \rp^2}
			      \rp                   \\
			       & = \frac \dh 4 \lp
			      1 - \frac 1 4 \sqrt{\abs{1 - \Fr^2}}
			      + 1 + \frac 1 4 \sqrt{\abs{1 - \Fr^2}}
			      \rp
			      = \frac {\dh} 2,
		      \end{aligned}
	      \end{equation*}
	      which completes the proof of \ref{item:def:H_property_WB}.

	\item Finally, we turn to \ref{item:def:properties_of_H_dry_wet}.
	      We first make the general remark that $\cH$, given by \eqref{eq:cH_version_2},
	      rewrites as
	      \begin{equation}
		      \label{eq:H_as_jump_times_bounded_function}
		      \cH (h_L, h_R, \bar q, \dZ) = (h_R - h_L) \mathcal{B}(h_L, h_R, \bar q, \dZ),
	      \end{equation}
	      with $\mathcal{B}$ a bounded function.
	      The boundedness of $\mathcal{B}$ comes from the boundedness assumption
	      on the Froude number, given in \eqref{eq:limh0}.

	      To prove the first equality of \ref{item:def:properties_of_H_dry_wet},
	      we consider the function
	      \begin{equation}
		      \label{eq:def_C_is_Fr_and_H}
		      \mathcal{C}(h_L, h_R, q_L) = 2 \Fr^2(h_L, h_R, q_L) \cH(h_L, h_R, q_L, \dZ).
	      \end{equation}
	      We have to prove that $\mathcal{C}$ goes to zero when $h_L$ tends to $0$,
	      in the two cases where $h_R = h_L$ and $h_R = 0$.
	      Technical Taylor expansions of $\mathcal{C}$ are provided in \cref{sec:Taylor_expansions_C}.
	      They show that this property is satisfied,
	      thereby proving the first equality of \ref{item:def:properties_of_H_dry_wet}.

	      To prove the second equality of \ref{item:def:properties_of_H_dry_wet},
	      we have to show that
	      \begin{equation}
		      \label{eq:H_limit_in_source_h_i_n}
		      \lim_{h_i^n \to 0^+} \frac {\cH \lp h_\imhp, h_\iphm, q_i, Z_\iph - Z_\imh \rp^3} {h_\imhp + h_\iphm} = 0.
	      \end{equation}
	      To that end, recall from property \ref{item:DW_req_1} that,
	      if $h_i = 0$, then $h_\imhp = h_\iphm = 0$.
	      This follows from the fact that~$\cH$ satisfies \ref{item:def:properties_of_H_dry_wet}.
	      Therefore, proving \eqref{eq:H_limit_in_source_h_i_n} requires proving that
	      \begin{equation}
		      \label{eq:H_limit_in_source}
		      \lim_{\substack{h_L \to 0^+ \\ h_R \to 0^+}}
		      \frac {\cH (h_L, h_R, \bar q, \dZ)^3} {h_L + h_R} = 0.
	      \end{equation}
	      Arguing \eqref{eq:H_as_jump_times_bounded_function}, we get
	      \begin{equation*}
		      \label{eq:H_limit_in_source_with_B}
                \left\lvert \frac {\cH^3} {h_L + h_R} \right\rvert
                \leq
                \frac {\left\lvert  h_R - h_L \right\rvert^3} {|h_L| + |h_R|}
                \left\lvert \mathcal{B} \right\rvert^3
                \leq
                \frac {\lp |h_R| + |h_L| \rp^3} {|h_L| + |h_R|}
                \max \left\lvert \mathcal{B} \right\rvert^3
                =
                \lp |h_R| + |h_L| \rp^2
                \max \left\lvert \mathcal{B} \right\rvert^3.
	      \end{equation*}
	      Taking the limit of the above expression when $h_L$ and $h_R$ go to zero
	      proves that \eqref{eq:H_limit_in_source} holds.
	      Therefore, the third equality of \ref{item:def:properties_of_H_dry_wet} is satisfied.

	      Finally, the proof of the third equality of \ref{item:def:properties_of_H_dry_wet}
	      is obtained by arguing the same arguments as above,
	      in addition to the well-known limit
	      \begin{equation*}
		      \lim_{\substack{x \to 0^+ \\ y \to 0^+}} \frac {x y} {x + y} = 0.
	      \end{equation*}

	      The proof of \ref{item:def:properties_of_H_dry_wet} is therefore concluded.
\end{enumerate}
We have therefore proven the three properties, which concludes the proof.
\end{proof}

As a conclusion, the expression \eqref{eq:cH_version_2} of $\cH$
satisfies the required properties, as proven in \cref{prop:expression_of_H}.
The only limit where this expression is ill-defined is when both
$\Fr^2$ tends to $1$ and $\dZ$ tends to $0$.
This corresponds to a well-known issue,
since the steady equations are themselves ill-defined when $\Fr^2 = 1$ and $\dZ = 0$.
This resonant case is well-documented in the literature,
as seen for instance in~\cite{CasCha2016,GomCasParRus2021} and references therein.
In the numerical experiments, we set $\cH = 0$ when $\Fr^2 = 1$ and $\dZ = 0$.

\section{\texorpdfstring{{\ra Low-cost high-order extension}}{Low-cost high-order extension}}
\label{sec:HO_WB}

In \cref{sec:scheme}, we proposed a {\ra hydrodynamic} reconstruction
technique able to turn any first-order,
non-well-balanced scheme into a fully well-balanced one.
The present section is dedicated to its high-order extension.
Since the first-order scheme {\ra does not require solving Bernoulli's relations,
		we aim to develop a high-order and well-balanced extension
		that does not require solving Bernoulli's relations either.}
Indeed, typically, fully well-balanced high-order extensions
require solving some nonlinear equations, which can be computationally expensive
and may fail in the case of multiple or non-existent solutions.
For instance, we mention \cite{CasGalLopPar2008,Xin2014,BriXin2020,GomCasParRus2021},
where nonlinear equations have to be numerically solved.

To avoid the computational expense associated with a nonlinear solver,
we adopt the method described in~\cite{BerBulFouMbaMic2022},
which was initially introduced on a generic hyperbolic system,
and later applied to the 2D shallow water case in~\cite{MicBerClaFou2021}.
For the sake of completeness,
the application of this technique to the current scheme is summarized below.
For the remainder of this section, we build a scheme of order $d + 1$.

To that end, we start by considering the following generic
polynomial reconstruction of degree $d$ in space
(see for instance \cite{DioClaLou2012,DioLouCla2013}):
\begin{equation}
	\label{eq:reconstruction}
	\widehat W_i^n(x) = W_i^n + \Pi_i^n (x - x_i),
\end{equation}
where the degree $d$ polynomial $\Pi_i^n$ is defined such that, for
all $x \in (x_\imh, x_\iph)$, the following relations hold:
\begin{equation}
	\label{eq:reconstruction_order}
	\widehat W_i^n(x) = W(x,t^n) + \cO(\dx^{d+1})
	\text{\quad and \quad}
	\frac 1 \dx \int_{x_\imh}^{x_\iph} \widehat W_i^n(x) \: dx = W_i^n.
\end{equation}
This reconstruction of degree $d$ naturally defines
a scheme of space order $d + 1$,
see for instance \cite{DioClaLou2012,DioLouCla2013}.
{\ra Note that this polynomial reconstruction involves
a usual slope limiter \cite{LeV2002}.
The limiters we use in practice are described in \cref{sec:numerics},
devoted to numerical simulations.}

Given the polynomial reconstruction \eqref{eq:reconstruction},
we can define a high-order, non-well-balanced scheme.
To that end, consider the following modification of \eqref{eq:scheme}:
\begin{equation}
	\label{eq:scheme_order_d}
	W_i^{n+1} =
	W_i^n
	- \frac \dt \dx \lp
	\cF \lp \widehat W_\iphm^n, \widehat W_\iphp^n \rp
	-
	\cF \lp \widehat W_\imhm^n, \widehat W_\imhp^n \rp
	\rp
	+ \dt \widehat \cS_i^n,
\end{equation}
where the polynomial reconstruction is applied at each inner interface, to get
\begin{equation*}
	\label{eq:reconstruction_high_order}
	\widehat W_\iphm^n
	=
	W_i^n
	+
	\Pi_i^n \lp \frac\dx 2 \rp
	\text{\quad and \quad}
	\widehat W_\iphp^n
	=
	W_\ip^n
	+
	\Pi_\ip^n \lp - \frac\dx 2 \rp,
\end{equation*}
and where the source term is nothing but a high-order approximation
of the averaged source term:
\begin{equation*}
	\label{eq:high_order_source_is_high_order_accurate}
	\widehat \cS_i^n
	=
	\frac 1 \dx \int_{x_\imh}^{x_\iph} S(W(x, t^n), x) \: dx + \cO(\dx^{d+1}).
\end{equation*}
To compute $\widehat \cS_i^n$ in practice,
one can use a quadrature formula of order $d + 1$, see for instance \cite{AbrSte1992}.

Equipped with the high-order scheme \eqref{eq:scheme_order_d}, we make it well-balanced,
without having to solve nonlinear equations.
We accomplish this by leveraging the procedure from \cite{BerBulFouMbaMic2022}.
To address this issue, we first introduce a steady solution detector.

\subsection{A steady state detector}

Let us define the following indicator,
which will help us write the well-balanced extension of the high-order scheme:
\begin{align}
	\label{eq:def_theta_iph}
	\theta_\iph^n      & =
	\frac{\varepsilon_\iph^n}{\varepsilon_\iph^n+\lp\dfrac\dx{C_\iph^n}\rp^{d+1}} \text{, with} \\
	\label{eq:def_epsilon}
	\varepsilon_\iph^n & =
	\left \lVert
	\begin{pmatrix}
		q_\ip^n \\
		\dfrac{(q_\ip^n)^2} {2 (h_\ip^n)^2} + g (h_\ip^n + Z_\ip)
	\end{pmatrix}
	-
	\begin{pmatrix}
		q_i^n \\
		\dfrac{(q_i^n)^2} {2 (h_i^n)^2} + g (h_i^n + Z_i)
	\end{pmatrix}
	\right \rVert,
\end{align}
where $C_{i+\frac 1 2}^n \neq 0$ is independent of $\dx$.
An expression of $C_{i+\frac 1 2}^n$ will be proposed
before performing the numerical experiments, in \eqref{eq:expression_of_C}.
The properties enjoyed by this indicator are summarized in the following result.

\begin{proposition}
	\label{prop:properties_theta}
	The expressions \eqref{eq:def_theta_iph} -- \eqref{eq:def_epsilon}
	ensure the following properties:
	\begin{enumerate}[label=(\thetheorem-\alph*)]
		\item if $(W_i^n, Z_i^n)$ and $(W_\ip^n, Z_\ip^n)$ define a steady state,
		      i.e. satisfy the local per interface steady state
		      relation~\eqref{eq:steadyperinterface},
		      then $\theta_\iph^n = 0$;
		\item otherwise, $\theta_\iph^n = 1 + \cO(\dx^{d+1})$.
	\end{enumerate}
\end{proposition}

\begin{proof}
	The proof is straightforward and present in \cite{BerBulFouMbaMic2022}.
	The first property is verified by inspection of \eqref{eq:def_epsilon},
	and the second one through a Taylor expansion of~$\theta_\iph^n$.
\end{proof}

\subsection{The high-order well-balanced reconstruction}

Equipped with the steady solution detector
\eqref{eq:def_theta_iph} -- \eqref{eq:def_epsilon}
and the high-order, non-well balanced scheme \eqref{eq:scheme_order_d},
we are able to fully state {\ra the high-order, well-balanced scheme}.

To that end, we modify the high-order scheme \eqref{eq:scheme_order_d} as follows:
\begin{equation}
	\label{eq:scheme_order_d_WB}
	W_i^{n+1} =
	W_i^n
	- \frac \dt \dx \lp
	\cF \lp \widetilde W_\iphm^n, \widetilde W_\iphp^n \rp
	-
	\cF \lp \widetilde W_\imhm^n, \widetilde W_\imhp^n \rp
	\rp
	+ \dt \widetilde \cS_i^n,
\end{equation}
where the polynomial reconstruction $\widehat W$
has been replaced with the modified polynomial reconstruction $\widetilde W$, defined by
\begin{equation*}
	\label{eq:reconstruction_high_order_WB}
	\widetilde W_\iphm^n
	=
	W_i^n
	+
	\theta_\iph^n \Pi_i^n \lp \frac\dx 2 \rp
	\text{\quad and \quad}
	\widetilde W_\iphp^n
	=
	W_\ip^n
	+
	\theta_\iph^n \Pi_\ip^n \lp - \frac\dx 2 \rp,
\end{equation*}
and where the high-order source term $\widehat \cS_i^n$ has been replaced with a
convex combination $\widetilde \cS_i^n$
between the high-order source term $\widehat \cS_i^n$
and the order one source term $\cS_i^n$, defined in \eqref{eq:source}:
\begin{equation*}
	\label{eq:source_high_order_WB}
	\widetilde \cS_i^n
	=
	\lp 1 - \frac {\theta_\imh^n + \theta_\iph^n} 2 \rp \cS_i^n
	+
	\frac {\theta_\imh^n + \theta_\iph^n} 2 \, \widehat \cS_i^n.
\end{equation*}
Note that \cref{prop:properties_theta} implies that,
if there is a steady state at the interface $x_{i+1/2}$,
then \smash{$\widetilde W_{i+1/2,-}^n = W_i^n$} and
\smash{$\widetilde W_{i+1/2,+}^n = W_\ip^n$}.
Otherwise, \smash{$\widetilde W_{i+1/2,-}^n$} and \smash{$\widetilde W_{i+1/2,+}^n$}
are high-order approximations of the solution at the interface $x_{i+1/2}$.

Equipped with this modified high-order reconstruction,
the high-order hydrodynamic reconstruction is simply computed by applying
\eqref{eq:hydrodynamic_reconstruction}
to \smash{$\widetilde W_{i+1/2,-}^n$} and
\smash{$\widetilde W_{i+1/2,+}^n$} at the interface $x_{i+1/2}$,
instead of applying it to $W_i$ and~$W_\ip$.

Thanks to these definitions, we are able to state the following result,
which represents a high-order extension of \cref{thm:scheme_properties}.

\begin{theorem}
	\label{thm:HO_scheme}
	The scheme \eqref{eq:scheme_order_d} enjoys the following properties:
	\begin{enumerate}[label=(\thetheorem-\alph*)]
		\item \label{item:thm:HO_scheme_properties_consistency}
		      it is high-order accurate,
		      i.e. consistent with the shallow water equations \eqref{eq:shallow_water}
		      up to $\cO(\dx^{d+1})$,
		\item \label{item:thm:HO_scheme_properties_positivity}
		      it is non-negativity-preserving,
		\item \label{item:thm:HO_scheme_properties_WB}
		      it is fully well-balanced, in the sense that it exactly preserves the steady states \eqref{eq:steady_discrete}.
	\end{enumerate}
\end{theorem}

\begin{proof}
	The proof of \cref{thm:HO_scheme} is present in \cite{BerBulFouMbaMic2022}.
	For the sake of conciseness, it is omitted here.
\end{proof}

\section{Numerical experiments}
\label{sec:numerics}

This last section comprises several numerical experiments,
designed to validate the properties of the scheme.
Firstly, in \cref{sec:order_of_accuracy}, the consistency and order of accuracy are assessed.
Secondly, in \cref{sec:numerical_well_balancedness}, we perform several experiments
to test the well-balanced property of the scheme
(for a lake at rest in \cref{sec:lake_at_rest}
and for moving steady states in \cref{sec:moving_steady_states}).
Thirdly, two dam-break problems are computed in \cref{sec:dam_breaks}.
Lastly, we simulate an unstable steady contact wave in \cref{sec:stationary_contact}.

Recall that any consistent and non-negativity-preserving
numerical flux $\cF$ can be used in the scheme \eqref{eq:scheme}.
For the numerical experiments, we use the simple HLL flux from \cite{HarLaxLee1983}.
This flux imposes a standard CFL restriction on the time step,
as discussed for instance in \cite{Vil1986,Tor2009}.
When it is required, in \cref{sec:order_of_accuracy} and in \cref{sec:dam_breaks},
the reference solution is given by the uncorrected HLL scheme
with a naive centered discretization of the source term.

For clarity, we label the schemes under consideration as follows:
\begin{itemize}[nosep]
	\item the HSR scheme is the first-order accurate hydrostatic reconstruction
	      from \cite{AudBouBriKlePer2004},
	\item the HDR$p$ scheme is the hydrodynamic reconstruction
	      \eqref{eq:hydrodynamic_reconstruction}
	      for the well-balanced scheme of order $p$,
	      constructed with the method described in \cref{sec:HO_WB} if $p \geq 2$.
	      For the second-order scheme, we use the \textsf{minmod} limiter
	      (see for instance \cite{Lee1979})
	      in conjunction with the SSPRK2 time discretization.
	      For the third-order scheme, we use the reconstruction from \cite{SchSeiTor2015}
	      with the SSPRK3 time discretization.
	      Both time discretizations are presented in,
	      for instance, \cite{GotShu1998,GotShuTad2001}.
\end{itemize}

Unless otherwise mentioned, the space domain is $(0, 1)$ and the gravity constant is $g = 9.81$.
As prescribed in~\cite{BerBulFouMbaMic2022},
the coefficient $C^n_\iph$ in \eqref{eq:def_theta_iph} is given by
$C_\iph^0 = 1$ and, for $n \geq 1$,
\begin{equation}
	\label{eq:expression_of_C}
	C_\iph^n
	=
	C_\theta \,
	\frac 1 2
	\left(
	\frac{\left\lVert W_{i+1}^n - W_{i+1}^{n-1} \right\rVert} {\Delta t}
	+
	\frac{\left\lVert W_i^n - W_i^{n-1} \right\rVert} {\Delta t}
	\right),
\end{equation}
where $C_\theta$ is a constant depending on the numerical experiment.
Unless mentioned in a specific experiment, we take $C_\theta = 1$.

Before proceeding with the numerical experiments themselves,
we briefly mention how dry zones are numerically treated.
Let $\varepsilon_m = 2^{-52}$ be the machine epsilon.
According to \eqref{eq:limh0}, the velocity $u$ is computed as follows:
\begin{equation*}
	u =
	\begin{dcases}
		\frac q h & \text{if } h > \varepsilon_m, \\
		0         & \text{otherwise.}
	\end{dcases}
\end{equation*}
In addition, for the high-order scheme,
the non-negativity-preserving limiting procedure from \cite{Ber2005}
is used on the reconstructed water height.
The remaining potential divisions by $h$ or $\dZ$ are handled by leveraging
properties \ref{item:def:H_property_dZ_0}, \ref{item:def:H_property_WB}
and \ref{item:def:properties_of_H_dry_wet}.

\begin{figure}[!htbp]

	\centering

	\begin{tikzpicture}

		\pgfmathsetmacro{\x}{2}
		\pgfmathsetmacro{\w}{2}


		\path[fill=fill_topo] (-\x,0) -- (0,0) -- (0,\x) -- (\x,\x) -- (\x,-\x/2) -- (-\x,-\x/2) -- cycle;
		\path[fill=white!50!blue, opacity=0.5]  (-\x,\x/2) -- (0,\x/2) -- (0,0) -- (-\x,0) -- cycle;

		\draw[thick] (-\x,0) -- (0,0) -- (0,\x) -- (\x,\x);
		\draw[thick] (-\x,\x/2) -- (0,\x/2);

		\draw[stealth-stealth] (-\x/2,\x/2) -- (-\x/2,0) node[left,midway]{$h_i$};
		\draw[stealth-stealth] (-\x/2,0) -- (-\x/2,-\x/2) node[left,midway]{$Z_i$};
		\draw[stealth-stealth] (\x/2,\x) -- (\x/2,-\x/2) node[right,midway]{$Z_\ip$};


		\path[fill=fill_topo] (\w+\x,\x) -- (\w+2*\x,\x) -- (\w+2*\x,0) -- (\w+3*\x,0) -- (\w+3*\x,-\x/2) -- (\w+\x,-\x/2) -- cycle;
		\draw[fill=white!50!blue, opacity=0.5]  (\w+2*\x,\x/2) -- (\w+3*\x,\x/2) -- (\w+3*\x,0) -- (\w+2*\x,0) -- cycle;

		\draw[thick] (\w+\x,\x) -- (\w+2*\x,\x) -- (\w+2*\x,0) -- (\w+3*\x,0);
		\draw[thick] (\w+2*\x,\x/2) -- (\w+3*\x,\x/2);

		\draw[stealth-stealth] (\w+3*\x-\x/2,\x/2) -- (\w+3*\x-\x/2,0) node[right,midway]{$h_\ip$};
		\draw[stealth-stealth] (\w+3*\x-\x/2,0) -- (\w+3*\x-\x/2,-\x/2) node[right,midway]{$Z_\ip$};
		\draw[stealth-stealth] (\w+\x+\x/2,\x) -- (\w+\x+\x/2,-\x/2) node[left,midway]{$Z_i$};

	\end{tikzpicture}

	\caption{%
		Non-smooth, emerged steady state at rest not governed by \eqref{eq:steady_with_constants}.
		Left panel: lake at rest with $h_\ip = 0$ and $h_i + Z_i < Z_\ip$.
		Right panel: lake at rest with $h_i = 0$ and $Z_i > h_\ip + Z_\ip$.
	}

	\label{fig:emerged_steady_state_drawing}

\end{figure}
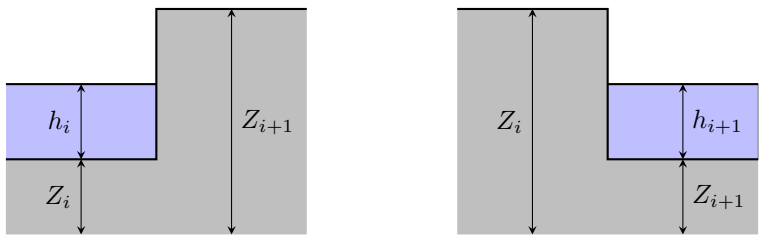

Further, to tackle steady solutions with an emerged bottom,
it should be noted that they
\emph{are not solution to Bernoulli's relation~\eqref{eq:steady_with_constants}}.
Indeed, after \cite{MicBerClaFou2016}, such steady states must be at rest,
which means that they must satisfy $q_0 = 0$ in \eqref{eq:steady_with_constants}.
In addition, the height and topography satisfy either
$h_i + Z_i < Z_\ip$ and $h_\ip = 0$ (\cref{fig:emerged_steady_state_drawing}, left panel), or
$Z_i > h_\ip + Z_\ip$ and $h_i = 0$ (\cref{fig:emerged_steady_state_drawing}, right panel).
Since such steady solutions are not solution to Bernoulli's relation,
$\cH$ cannot capture them without a modification.
To address this issue, we impose the following additional properties on $\cH$:
\begin{equation*}
	\label{eq:prop:properties_of_H_dry_wet_emerged}
	\begin{dcases}
		h_L < \dZ \quad   & \implies \quad \cH(h_L, 0, 0, \dZ) = \frac {\dh} 2,  \\
		h_R < - \dZ \quad & \implies \quad \cH(0, h_R, 0, \dZ) = \frac {\dh} 2.
	\end{dcases}
\end{equation*}
These properties ensure that steady states at rest
with an emerged bottom are correctly captured by the scheme.
Note that a similar technique was already used in
\cite{MicBerClaFou2016,MicBerClaFou2017,MicBerClaFou2021}
to make smooth-steady-state-capturing schemes able to
handle such non-smooth steady solutions at rest.

\subsection{Order of accuracy}
\label{sec:order_of_accuracy}

The first round of numerical experiments consists in measuring the order of accuracy.
To that end, we introduce this useful compactly supported~$\mathcal{C}^\infty$ bump function:
\begin{equation*}
	\omega(x) =
	\begin{dcases}
		\exp \left( 1 - \frac 1 {1 - (4 (x - 1/2))^2} \right) & \text{if } \lvert x - 1/2 \rvert < 1/4, \\
		0                                                     & \text{otherwise.}
	\end{dcases}
\end{equation*}
Note that $\omega(x)$ vanishes when $\lvert x - 1/2 \rvert \geq 1/4$.
We take $Z(x) = \omega(x)$, and
the initial condition is given by
\begin{equation*}
	h(x, 0) = 2 - Z(x) + \cos^2(2 \pi x)
	\text{ \quad and \quad }
	q(x, 0) = \sin(2 \pi x).
\end{equation*}
These initial data are evolved until the final time $t_{\text{end}} = 5 \cdot 10^{-3}$,
chosen in order to avoid the formation of discontinuities,
so as to be able to compute the order of accuracy.
A reference solution, to which the results of the schemes are compared,
is computed with $20 \cdot 2^{12} = 81920$ cells.
Periodic boundary conditions are prescribed.

\begin{figure}[!ht]
	\centering
	\includegraphics{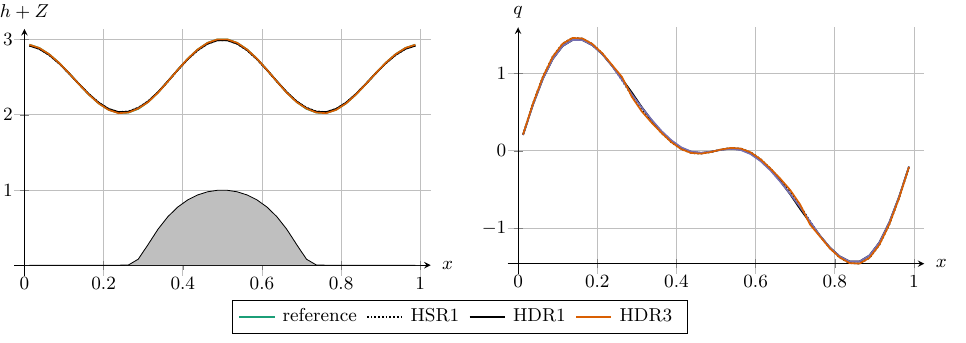}
	\caption{%
		Experiment from \cref{sec:order_of_accuracy}:
		values of $h$ (left panel) and $q$ (right panel)
		at time $t_{\text{end}}$ with $40$ cells.
	}
	\label{fig:SWT_dyadic_hydro_solution}
\end{figure}

In \cref{fig:SWT_dyadic_hydro_solution}, we display the reference solution and
the approximations given by the HSR$1$, HDR$1$ and HDR$3$ schemes with $40$ cells.
We observe good agreement with the exact solution.

\begin{figure}[!ht]
	\centering
	\includegraphics{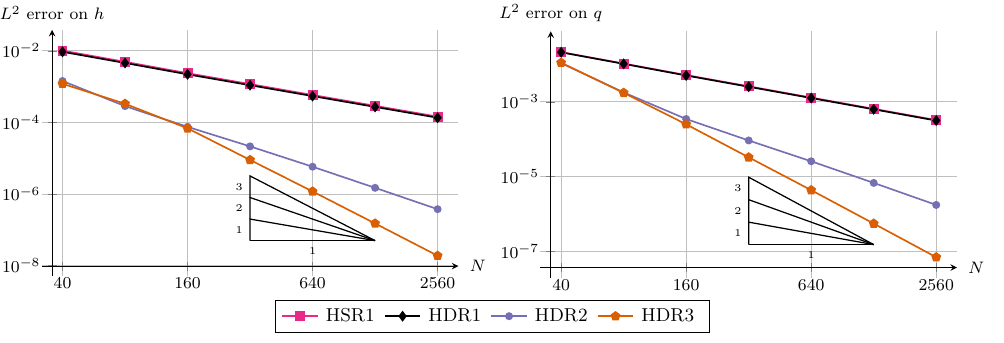}
	\caption{%
		Experiment from \cref{sec:order_of_accuracy}:
		error in $L^2$ norm
		on $h$ (left panel)
		and on $q$ (right panel)
		with respect to the number of cells.
	}
	\label{fig:SWT_dyadic_hydro_errors}
\end{figure}

To obtain a more precise assessment of the error,
error lines are shown in \cref{fig:SWT_dyadic_hydro_errors}.
We observe that the schemes exhibit the expected orders of accuracy.
As expected, neither the high-order well-balanced procedure
nor the hydrodynamic reconstruction impedes the order of accuracy.
For the sake of completeness,
the values of the errors are reported in \cref{tab:SWT_dyadic_errors}
(we only present errors on $h$, but the results for $q$ are similar).

\begin{table}[!ht]
	\newcommand\myCaption{%
		Experiment from \cref{sec:order_of_accuracy}:
		errors (in $L^2$ norm) and orders of accuracy on $h$.
		The errors on $q$ are displayed in \cref{fig:SWT_dyadic_hydro_errors}.}
	\centering
	\begin{filecontents*}{SWT_dyadic_hydro_errors.csv}
%% LaTeX2e file `SWT_dyadic_hydro_errors.csv'
%% generated by the `filecontents' environment
%% from source `Berthon_Michel-Dansac' on 2025/01/07.
%%
    N h_HSR_o1_error h_HSR_o1_order q_HSR_o1_error q_HSR_o1_order h_HDR_o1_error h_HDR_o1_order q_HDR_o1_error q_HDR_o1_order h_HDR_o2theta_error h_HDR_o2theta_order q_HDR_o2theta_error q_HDR_o2theta_order h_HDR_o3theta_error h_HDR_o3theta_order q_HDR_o3theta_error q_HDR_o3theta_order
    2.000000000000000000e+01 2.145575615791598181e-02 nan 4.352007220791001552e-02 nan 1.979453392090003289e-02 nan 4.347030202347802047e-02 nan 5.000447967507661538e-03 nan 3.238599573192246089e-02 nan 3.128540253803204548e-03 nan 2.657036370378268800e-02 nan
    4.000000000000000000e+01 1.010141914244681340e-02 1.086806755797434354e+00 2.151854599468631946e-02 1.016100351815892200e+00 9.360541468362732279e-03 1.080438207050477439e+00 2.129949606733383752e-02 1.029210821625882177e+00 1.435504905287514647e-03 1.800499084787247961e+00 1.094958837280375109e-02 1.564493465812479833e+00 1.197121368309377347e-03 1.385920241680359322e+00 1.119656147755750120e-02 1.246762234914448486e+00
    8.000000000000000000e+01 4.904948787631045574e-03 1.042248012079331598e+00 1.066052028503623739e-02 1.013302748151616139e+00 4.555841960300561003e-03 1.038874284524214486e+00 1.044121015708946928e-02 1.028530364766286009e+00 2.848928801582346752e-04 2.333066788001756819e+00 1.790523304400369522e-03 2.612423434535578348e+00 3.318169279942940076e-04 1.851110031205681627e+00 1.749794786004591320e-03 2.677798101714361945e+00
    1.600000000000000000e+02 2.365421719264307332e-03 1.052140655989880491e+00 5.231033626207548005e-03 1.027109901210752296e+00 2.197880450253883369e-03 1.051604785976995027e+00 5.094796138868472671e-03 1.035192607595032177e+00 7.508514266783966325e-05 1.923820196786918180e+00 3.478323428331128774e-04 2.363917304541705011e+00 6.838366436840719758e-05 2.278663851399910012e+00 2.508394861216184249e-04 2.802349360044251281e+00
    3.200000000000000000e+02 1.173902474439181840e-03 1.010784860317662126e+00 2.616460810829676544e-03 9.994793934720521689e-01 1.092395575592094515e-03 1.008617539583577916e+00 2.542608085074353399e-03 1.002715316450360605e+00 2.138798919278802165e-05 1.811726612531691094e+00 9.248315466234220726e-05 1.911129571313552544e+00 8.955065155612604666e-06 2.932875898105152146e+00 3.269156242103253248e-05 2.939774234855239587e+00
    6.400000000000000000e+02 5.829112545679652081e-04 1.009964395213373445e+00 1.305360179435241389e-03 1.003168715464206606e+00 5.427120680948985254e-04 1.009236481195676660e+00 1.267016592077034034e-03 1.004873686416923206e+00 5.790366690758137420e-06 1.885074234047782626e+00 2.562836558848393388e-05 1.851449134174108613e+00 1.180255794051072802e-06 2.923604363702692055e+00 4.316343379789226752e-06 2.921036783377485069e+00
    1.280000000000000000e+03 2.903908598120166860e-04 1.005280212222042779e+00 6.518806276220300575e-04 1.001768227618294649e+00 2.704312406786839730e-04 1.004925165890679040e+00 6.323607221538457253e-04 1.002615751838804803e+00 1.490025180764532356e-06 1.958318001904883587e+00 6.750150929810501639e-06 1.924749810261773275e+00 1.506924163981087250e-07 2.969420845174138179e+00 5.488442359971377666e-07 2.975340970386263884e+00
    2.560000000000000000e+03 1.449256533465902207e-04 1.002683055035057391e+00 3.257409526433891226e-04 1.000882693401382806e+00 1.349809054133168696e-04 1.002506487779823097e+00 3.158947746854309790e-04 1.001303687095017825e+00 3.775555460379693010e-07 1.980575897087887727e+00 1.733961178889268586e-06 1.960848161547015689e+00 1.897287164061326377e-08 2.989596855456701174e+00 6.885326240413867458e-08 2.994799846342502736e+00
    5.120000000000000000e+03 7.240279124760002206e-05 1.001195767734669673e+00 1.628322111187340015e-04 1.000338990256651739e+00 6.743926717296166767e-05 1.001094571592347826e+00 1.578875635172047971e-04 1.000546535696498873e+00 9.502396526455631158e-08 1.990325594115825325e+00 4.395316413423000853e-07 1.980032666814935949e+00 2.378406239281197481e-09 2.995870996404891518e+00 8.612953179799277407e-09 2.998945119853682950e+00
    1.024000000000000000e+04 3.618025520657260633e-05 1.000842732256411516e+00 8.139728116702476884e-05 1.000333607182393836e+00 3.370062542347098007e-05 1.000813493856001601e+00 7.891965649809868917e-05 1.000440955731664694e+00 2.385757168227498366e-08 1.993844203055196385e+00 1.107540229846404118e-07 1.988607923341653816e+00 2.977112548149057810e-10 2.998009484955318271e+00 1.076882068126200399e-09 2.999647723141365407e+00
\end{filecontents*}
\pgfplotstabletypeset[
		columns={N, h_HSR_o1_error, h_HSR_o1_order, h_HDR_o1_error, h_HDR_o1_order, h_HDR_o2theta_error, h_HDR_o2theta_order, h_HDR_o3theta_error, h_HDR_o3theta_order},
		skip rows between index={0}{1},
		skip rows between index={8}{10},c
		every row 1 column h_HSR_o1_order/.style={postproc cell content/.style={@cell content=---}},
		every row 1 column h_HDR_o1_order/.style={postproc cell content/.style={@cell content=---}},
		every row 1 column h_HDR_o2theta_order/.style={postproc cell content/.style={@cell content=---}},
		every row 1 column h_HDR_o3theta_error/.style={postproc cell content/.style={@cell content=---}},
		every head row/.append style={before row={%
                        \toprule
						\multicolumn{1}{c}{} &
						\multicolumn{2}{c}{$h$, HSR, $\mathbb{P}_0$} &
						\multicolumn{2}{c}{$h$, HDR, $\mathbb{P}_0$} &
						\multicolumn{2}{c}{$h$, HDR, $\mathbb{P}_1$} &
						\multicolumn{2}{c}{$h$, HDR, $\mathbb{P}_2$} \\
						\cmidrule(lr){2-3}
						\cmidrule(lr){4-5}
						\cmidrule(lr){6-7}
						\cmidrule(lr){8-9}
					}}
	]{SWT_dyadic_hydro_errors.csv}
	\caption{\myCaption}
	\label{tab:SWT_dyadic_errors}
\end{table}

\subsection{Well-balanced property}
\label{sec:numerical_well_balancedness}

We now turn to experiments that assess the well-balanced property.
We first examine submerged and emerged lake at rest steady solutions
in \cref{sec:lake_at_rest}.
Next, in \cref{sec:moving_steady_states}, we tackle the experiments from~\cite{GouMau1997},
which involve moving steady solutions that are reached after a transient, unsteady state.

\subsubsection{Lake at rest}
\label{sec:lake_at_rest}

We begin by studying steady states at rest, taking $Z(x) = \omega(x)$ once again.
The initial discharge is zero everywhere ($q(x, 0) = 0$),
and the initial height is given in \cref{tab:lake_at_rest_setup}.
Note that the resulting initial condition is nothing but
a steady state at rest of the shallow water system \eqref{eq:shallow_water}.
Therefore, since the HSR and HDR schemes are well-balanced,
we expect them to exactly preserve this initial condition.
We fix the final time at $t_{\text{end}} = 1$,
take $50$ cells,
and prescribe the exact steady solution as inhomogeneous Dirichlet boundary conditions.
We conduct two experiments:
the first one has a submerged bottom (no dry zones),
and the second one has an emerged bottom (with a dry area).

\begin{table}[!ht]
	\newcommand\myCaption{%
		Setup of the experiments from \cref{sec:lake_at_rest}.
	}
	\centering
	\begin{tabular}{rcccc}
		\toprule
		experiment       & figure                                 & $h(x, 0)$             \\
		\midrule
		submerged bottom & \cref{fig:SWT_rest_hydro_solution}     & $2 - Z(x)$            \\
		emerged bottom   & \cref{fig:SWT_rest_dry_hydro_solution} & $\max(0, 0.5 - Z(x))$ \\
		\bottomrule
	\end{tabular}
	\caption{\myCaption}
	\label{tab:lake_at_rest_setup}
\end{table}

\begin{figure}[!ht]
	\centering
	\includegraphics{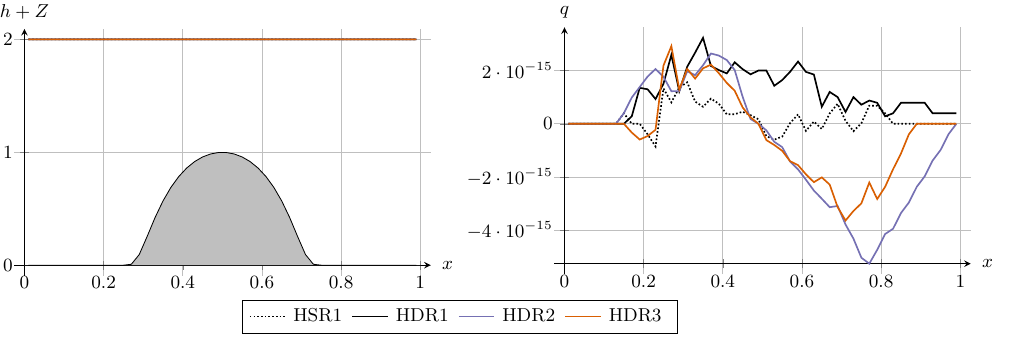}
	\caption{%
		Lake at rest with submerged bottom from \cref{sec:lake_at_rest}:
		free surface $h+Z$ (left panel)
		and discharge~$q$ (right panel),
		displayed at time $t_{\text{end}} = 1$ with $50$ cells.
	}
	\label{fig:SWT_rest_hydro_solution}
\end{figure}

\paragraph*{Submerged bottom.}
First, in \cref{fig:SWT_rest_hydro_solution},
we display the results of the lake at rest with a submerged bottom.
As expected, the initial condition is exactly preserved (up to machine precision)
by all the schemes under consideration (HSR scheme, HDR scheme, and its high-order extensions).
These conclusions are confirmed by the values of the errors
reported in \cref{tab:SWT_rest_hydro_errors}.

\begin{table}[!ht]
	\newcommand\myCaption{%
		Lake at rest with submerged bottom from \cref{sec:lake_at_rest}:
		errors (in $L^2$ norm) between the initial condition and the approximate solutions
		at time $t_{\text{end}} = 1$,
		using $50$ cells.
	}
	\centering
	\begin{filecontents*}{SWT_rest_hydro_errors.csv}
%% LaTeX2e file `SWT_rest_hydro_errors.csv'
%% generated by the `filecontents' environment
%% from source `Berthon_Michel-Dansac' on 2025/01/07.
%%
    variable HSR_o1 HDR_o1 HDR_o2theta HDR_o3theta
    h 8.881784197001253e-17 2.010699415441723e-16 1.0877919644084146e-16 4.4408920985006264e-17
    q 5.245806859157568e-16 1.4170880853618084e-15 2.3201106970426504e-15 1.605152245019692e-15
\end{filecontents*}
\pgfplotstabletypeset[
		columns={variable, HSR_o1, HDR_o1, HDR_o2theta, HDR_o3theta},
	]{SWT_rest_hydro_errors.csv}
	\caption{\myCaption}
	\label{tab:SWT_rest_hydro_errors}
\end{table}

\begin{figure}[!ht]
	\centering
	\includegraphics{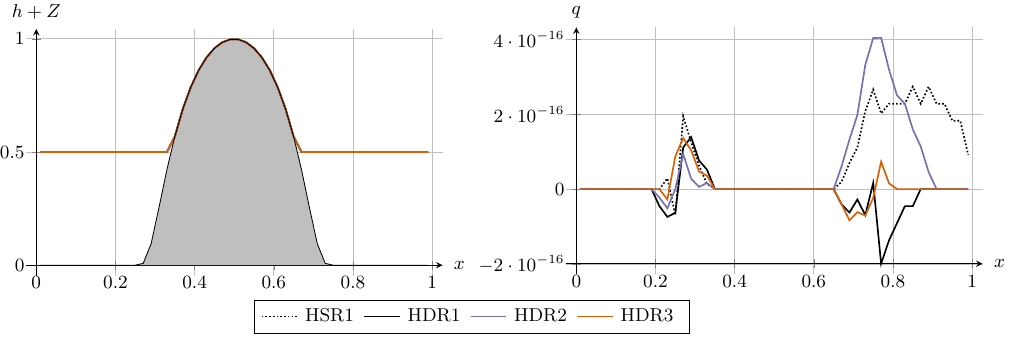}
	\caption{%
		Lake at rest with emerged bottom from \cref{sec:lake_at_rest}:
		free surface $h+Z$ (left panel)
		and discharge $q$ (right panel),
		displayed at time $t_{\text{end}} = 1$ with $50$ cells.
	}
	\label{fig:SWT_rest_dry_hydro_solution}
\end{figure}

\paragraph*{Emerged bottom.}
Next, the results of the steady state at rest with emerged bottom
are depicted in \cref{fig:SWT_rest_dry_hydro_solution},
and the errors are collected in \cref{tab:SWT_rest_dry_hydro_errors}.
Similarly to the submerged bottom case,
the dry zones did not negatively impact the well-balanced property.

\begin{table}[!ht]
	\newcommand\myCaption{%
		Lake at rest with emerged bottom from \cref{sec:lake_at_rest}:
		errors (in $L^2$ norm) between the initial condition and the approximate solutions
		at time $t_{\text{end}} = 1$,
		using $50$ cells.
	}
	\centering
	\begin{filecontents*}{SWT_rest_dry_hydro_errors.csv}
%% LaTeX2e file `SWT_rest_dry_hydro_errors.csv'
%% generated by the `filecontents' environment
%% from source `Berthon_Michel-Dansac' on 2025/01/07.
%%
    variable HSR_o1 HDR_o1 HDR_o2theta HDR_o3theta
h 1.851527828850887e-17 2.754662226278516e-17 3.065703529144609e-17 1.3165625507162628e-17
q 1.2405118483756166e-16 5.172192065851505e-17 1.2363777419240397e-16 3.593435675428176e-17
\end{filecontents*}
\pgfplotstabletypeset[
		columns={variable, HSR_o1, HDR_o1, HDR_o2theta, HDR_o3theta},
	]{SWT_rest_dry_hydro_errors.csv}
	\caption{\myCaption}
	\label{tab:SWT_rest_dry_hydro_errors}
\end{table}

\subsubsection{Moving steady solutions}
\label{sec:moving_steady_states}

To assess the ability of the HDR scheme to capture moving steady solutions,
we now examine the well-known test cases from \cite{GouMau1997}.
Namely, we run three test cases:
a subcritical flow, a transcritical flow without shock and a transcritical flow with a shock.
Each of these test cases follows the same principle:
the initial condition consists in a steady state at rest,
which is then perturbed by an inflow boundary condition at the left of the domain.
After a transient state,
the resulting flow becomes a moving steady state (with nonzero velocity).
For the subcritical flow and the transcritical flow without shock,
this moving steady state satisfies
\begin{equation*}
	q = \cst = q_0
	\text{\qquad and \qquad}
	B = \frac{q^2}{2 h^2} + g (h + Z) = \cst = \frac{q_0^2}{2 H_0^2} + g H_0.
\end{equation*}
This final steady state therefore depends on two parameters:
the inflow discharge, denoted by $q_0$, and the initial free surface, denoted by $H_0$.
We expect the HDR scheme and its high-order extensions
to exactly capture the final steady state,
and the HSR scheme to have a nonzero approximation error.
However, the transcritical flow with shock, being a non-smooth steady state,
is not expected to be exactly captured by the HDR scheme.

The space domain in $(0, 25)$,
where the function $Z(x) = \max(0, 0.05 (x - 8) (12 - x))$ is considered.
The initial conditions are defined as $h(x, 0) + Z(x) = H_0$ and $q(x, 0) = q_0$,
with $H_0$ and $q_0$ provided in \cref{tab:moving_steady_states_setup}.
The final times are also provided in \cref{tab:moving_steady_states_setup},
and we take $75$ discretization cells.
At the left boundary, we prescribe homogeneous Neumann boundary conditions on $h$,
and we impose $q(0, t) = q_0$.
At the right boundary, we prescribe homogeneous Neumann boundary conditions on $q$,
and we impose $h(25, t) = H_0$ if the flow is subcritical;
otherwise, homogeneous Neumann boundary conditions are prescribed on $h$.

\begin{table}[!ht]
	\newcommand\myCaption{%
		Setup of the experiments from \cref{sec:moving_steady_states}.
	}
	\centering
	\begin{tabular}{rcccc}
		\toprule
		experiment               & figure                            & $q_0$  & $H_0$  & $t_{\text{end}}$ \\
		\midrule
		subcritical              & \cref{fig:SWT_GM3_hydro_solution} & $4.42$ & $2$    & \;$500$\;        \\
		transcritical            & \cref{fig:SWT_GM2_hydro_solution} & $1.53$ & $0.66$ & \;$125$\;        \\
		transcritical with shock & \cref{fig:SWT_GM1_hydro_solution} & $0.18$ & $0.33$ & \;$1000$\;       \\
		\bottomrule
	\end{tabular}
	\caption{\myCaption}
	\label{tab:moving_steady_states_setup}
\end{table}

{\ra In this section, the exact solution satisfies $q = \cst$ and $B = \cst$.
As a consequence, the errors are evaluated according to $q$ and $B$.
Namely, with $N$ the number of cells in the mesh, we compute
\begin{equation*}
	e_q = \sqrt{\frac 1 \dx \sum_{i=1}^{N-1} | q_{i+1}^n - q_i^n |^2}
	\text{\qquad and \qquad}
	e_B = \sqrt{\frac 1 \dx \sum_{i=1}^{N-1} | B_{i+1}^n - B_i^n |^2}.
\end{equation*}
}

\begin{figure}[!ht]
	\centering
	\includegraphics{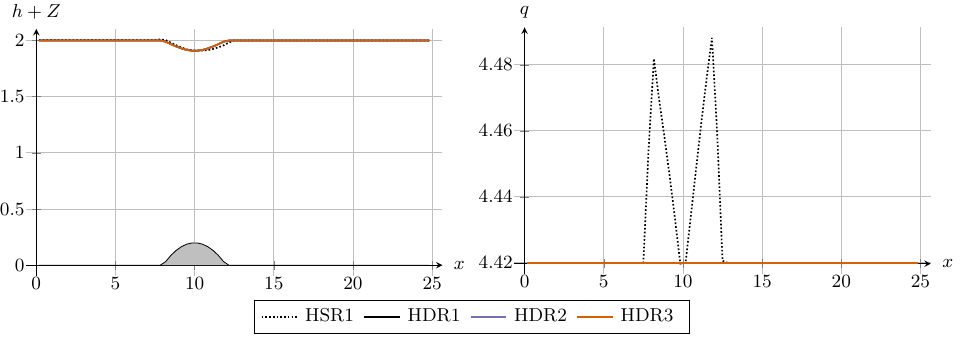}
	\caption{%
		Subcritical flow from \cref{sec:moving_steady_states}:
		free surface $h+Z$ (left panel)
		and discharge $q$ (right panel),
		displayed at time $t_{\text{end}} = 500$ with $75$ cells.
	}
	\label{fig:SWT_GM3_hydro_solution}
\end{figure}

\paragraph*{Subcritical flow.}
The first experiment, which converges towards a subcritical flow,
is illustrated in \cref{fig:SWT_GM3_hydro_solution}.
\Cref{tab:SWT_GM3_hydro_errors} contains the values
of the errors to the underlying steady state.
As expected, we observe that the HDR$1$ scheme, contrary to the HSR$1$ scheme,
exactly captures the resulting moving steady state.
In addition, the high-order extensions HDR$2$ and HDR$3$
also exactly capture the moving steady state.

\begin{table}[!ht]
	\newcommand\myCaption{%
		Subcritical flow from \cref{sec:moving_steady_states}:
		errors (in $L^2$ norm) between the exact steady state and the approximate solutions
		at time $t_{\text{end}} = 500$,
		using $75$ cells.
	}
	\centering
	\begin{filecontents*}{SWT_GM3_hydro_errors.csv}
%% LaTeX2e file `SWT_GM3_hydro_errors.csv'
%% generated by the `filecontents' environment
%% from source `Berthon_Michel-Dansac' on 2025/01/07.
%%
    variable HSR_o1 HDR_o1 HDR_o2theta HDR_o3theta
q 0.07726987694555672 1.0583867371683362e-14 1.3143840315653625e-14 1.2962537282957974e-14
B 0.1789217705754612 2.7288911567313107e-14 3.605610588796615e-14 2.6822400070865448e-14
\end{filecontents*}
\pgfplotstabletypeset[
		columns={variable, HSR_o1, HDR_o1, HDR_o2theta, HDR_o3theta},
	]{SWT_GM3_hydro_errors.csv}
	\caption{\myCaption}
	\label{tab:SWT_GM3_hydro_errors}
\end{table}

\begin{figure}[!ht]
	\centering
	\includegraphics{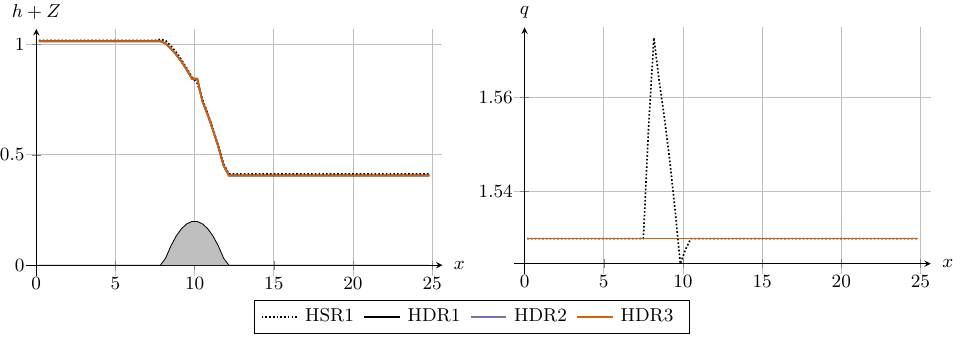}
	\caption{%
		Transcritical flow without shock from \cref{sec:moving_steady_states}:
		free surface $h+Z$ (left panel)
		and discharge $q$ (right panel),
		displayed at time $t_{\text{end}} = 125$ with $75$ cells.
	}
	\label{fig:SWT_GM2_hydro_solution}
\end{figure}

\paragraph*{Transcritical flow.}
The results of the second experiment,
which involves a transcritical steady flow,
are displayed in \cref{fig:SWT_GM2_hydro_solution}
and \cref{tab:SWT_GM2_hydro_errors}.
Similar to the previous case,
we observe that the steady state is exactly captured by the HDR scheme,
unlike the HSR scheme.
However, we observe a small kink near $x = 10$.
This defect arises because, at the critical point $x = 10$,
the Froude number is $1$ and the topography derivative vanishes.
Note that similar defects were already observed in earlier work,
see for instance~\cite{MicBerClaFou2016,BerMbaLeSec2021}.
{\ra This amplitude of the kink is reduced with the mesh size.
It is worth noting that,
although the numerical solution of the HDR scheme presents this kink,
it still satisfies $q = \cst$ and $B = \cst$ over the whole space domain,
as evidenced by \cref{tab:SWT_GM3_hydro_errors}.
}


\begin{table}[!ht]
	\newcommand\myCaption{%
		Transcritical flow without shock from \cref{sec:moving_steady_states}:
		errors (in $L^2$ norm) between the exact steady state and the approximate solutions
		at time $t_{\text{end}} = 125$,
		using $75$ cells.
	}
	\centering
	\begin{filecontents*}{SWT_GM2_hydro_errors.csv}
%% LaTeX2e file `SWT_GM2_hydro_errors.csv'
%% generated by the `filecontents' environment
%% from source `Berthon_Michel-Dansac' on 2025/01/07.
%%
    variable HSR_o1 HDR_o1 HDR_o2theta HDR_o3theta
q 0.037423772530638985 4.729393712878033e-14 5.1465992662138324e-14 5.214787393099677e-14
B 0.14502577105849515 4.504387026253897e-14 5.120716403936702e-14 5.91822812941618e-14
\end{filecontents*}
\pgfplotstabletypeset[
		columns={variable, HSR_o1, HDR_o1, HDR_o2theta, HDR_o3theta},
	]{SWT_GM2_hydro_errors.csv}
	\caption{\myCaption}
	\label{tab:SWT_GM2_hydro_errors}
\end{table}

\begin{figure}[!ht]
	\centering
	\includegraphics{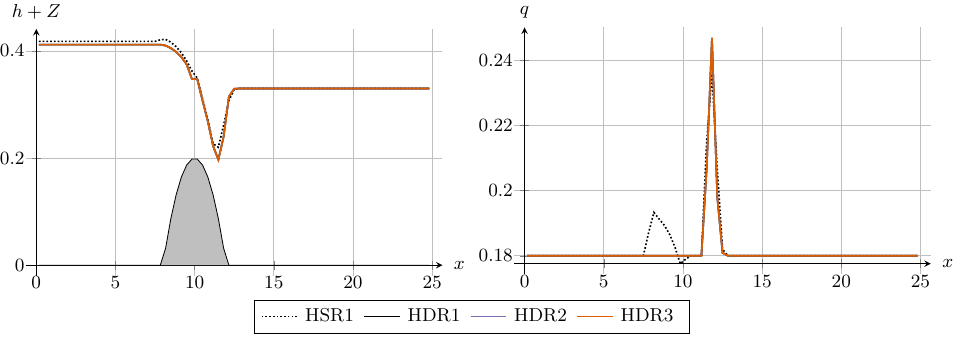}
	\caption{%
		Transcritical flow with shock from \cref{sec:moving_steady_states}:
		free surface $h+Z$ (left panel)
		and discharge $q$ (right panel),
		displayed at time $t_{\text{end}} = 1000$ with $75$ cells.
	}
	\label{fig:SWT_GM1_hydro_solution}
\end{figure}

\paragraph*{Transcritical flow with shock.}
The results of the final experiment are displayed in \cref{fig:SWT_GM1_hydro_solution}.
This experiment involves a transcritical flow with a shock;
as expected, since it is not smooth, it is not exactly captured by the HDR scheme,
let alone by the HSR scheme.
However, note that the loss of precision of the HDR scheme
only occurs in the vicinity of the shock (around $x = 12$),
while the continuous steady states before and after the shock are exactly preserved.

\begin{table}[!ht]
	\newcommand\myCaption{%
		Transcritical flow with shock from \cref{sec:moving_steady_states}:
		errors (in $L^2$ norm) between the exact steady state and the approximate solutions
		at time $t_{\text{end}} = 1000$,
		using $75$ cells.
	}
	\centering
	\begin{filecontents*}{SWT_GM1_hydro_errors.csv}
%% LaTeX2e file `SWT_GM1_hydro_errors.csv'
%% generated by the `filecontents' environment
%% from source `Berthon_Michel-Dansac' on 2025/01/07.
%%
    variable HSR_o1 HDR_o1 HDR_o2theta HDR_o3theta
q 0.05723156534824204 0.07085475615353774 0.07085550790307049 0.0708550986065084
B 0.692215945309875 0.8167153662406544 0.8167297010393522 0.8167218964027515
\end{filecontents*}
\pgfplotstabletypeset[
		columns={variable, HSR_o1, HDR_o1, HDR_o2theta, HDR_o3theta},
	]{SWT_GM1_hydro_errors.csv}
	\caption{\myCaption}
	\label{tab:SWT_GM1_hydro_errors}
\end{table}

\subsection{Dam-break problems}
\label{sec:dam_breaks}

The purpose of these experiments is to evaluate the performance of the scheme
on two standard dam-break problems:
the first one without dry areas, and the second one with a dry area.
For both problems, we take $Z(x) = x / 2$ and set the initial discharge to $q(x, 0) = 0$.
The initial water height is determined according to \cref{tab:dam_breaks_setup},
which also contains the values of $t_{\text{end}}$ and $\theta_C$.
Homogeneous Neumann boundary conditions are prescribed,
and we take~$50$ discretization cells.

\begin{table}[!ht]
	\newcommand\myCaption{%
		Setup of the experiments from \cref{sec:dam_breaks}.
	}
	\centering
	\begin{tabular}{rcccc}
		\toprule
		experiment
		 & figure
		 & $h(x, 0) + Z(x)$
		 & $t_{\text{end}}$                                                                & $\theta_C$ \\
		\midrule
		wet dam-break
		 & \cref{fig:SWT_RP_wet_hydro_solution}
		 & $\begin{dcases} 1.5 & \text{if } x < 0.5 \\ 1 & \text{otherwise} \end{dcases}$
		 & \;$0.05$\;
		 & \;$0.15$\;                                                                                   \\
		dry dam-break
		 & \cref{fig:SWT_RP_dry_hydro_solution}
		 & $\begin{dcases} 1 & \text{if } x < 0.5 \\ Z(x) & \text{otherwise} \end{dcases}$
		 & \;$0.075$\;
		 & \;$0.1$\;                                                                                    \\
		\bottomrule
	\end{tabular}
	\caption{\myCaption}
	\label{tab:dam_breaks_setup}
\end{table}

\begin{figure}[!ht]
	\centering
	\includegraphics{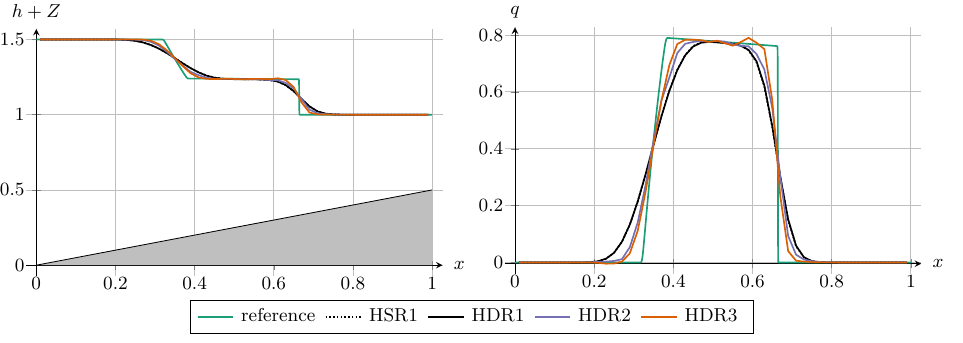}
	\caption{%
		Wet dam-break from \cref{sec:dam_breaks}:
		free surface $h+Z$ (left panel)
		and discharge $q$ (right panel),
		displayed at time $t_{\text{end}} = 0.05$ with $50$ cells.
	}
	\label{fig:SWT_RP_wet_hydro_solution}
\end{figure}

\paragraph*{Wet dam-break.}
\Cref{fig:SWT_RP_wet_hydro_solution} depicts the solutions of the wet dam-break experiment.
We observe no difference between the HSR$1$ and HDR$1$ schemes,
and it is worth noting that the HDR scheme's high-order extensions
provide a more accurate approximation of the exact solution,
despite minor oscillations on the discharge for the HDR$3$ scheme.
These oscillations are solely due to the high-order polynomial reconstruction.

\begin{figure}[!ht]
	\centering
	\includegraphics{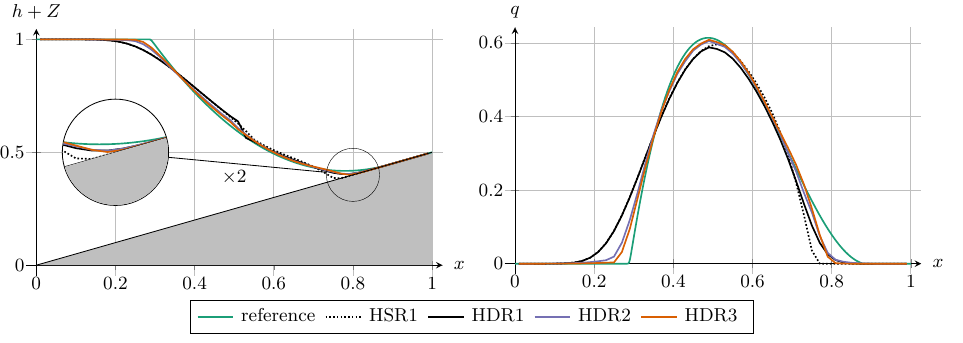}
	\caption{%
		Dry dam-break from \cref{sec:dam_breaks}:
		free surface $h+Z$ (left panel)
		and discharge $q$ (right panel),
		displayed at time $t_{\text{end}} = 0.075$ with $50$ cells.
	}
	\label{fig:SWT_RP_dry_hydro_solution}
\end{figure}

\paragraph*{Dry dam-break.}
The results of the dry dam-break problem are presented in \cref{fig:SWT_RP_dry_hydro_solution}.
Two significant differences between the HSR$1$ and HDR$1$ schemes are noteworthy.
First, the HDR$1$ scheme produces a small kink near the critical point $x = 0.5$.
However, that this kink disappears when the mesh is refined,
or when increasing the order of accuracy by using the HDR$2$ or HDR$3$ schemes.
Second, the HDR scheme produces a more accurate approximation of the wet/dry transition,
than the HSR scheme, despite having the same number of cells.

\subsection{Stationary contact wave}
\label{sec:stationary_contact}

This last experiment corresponds to the situation described in
\cref{sec:comments_steady}.
The discontinuous topography is given by
\begin{equation*}
	\label{eq:topography_stationary_contact}
	Z(x)
	=
	\begin{dcases}
		Z_L = 0    & \text{if } x < 0.5, \\
		Z_R = 0.01 & \text{otherwise.}
	\end{dcases}
\end{equation*}
We consider a Riemann problem with initial data
\begin{equation*}
	\label{eq:initial_data_stationary_contact}
	W(x, 0)
	=
	\begin{dcases}
		(h_L, q_L)^\intercal & \text{if } x < 0.5, \\
		(h_R, q_R)^\intercal & \text{otherwise,}
	\end{dcases}
\end{equation*}
with $q_L = q_R = 1$, $h_L = 1$ and
where $h_R \simeq 0.2545853624828563$ is defined such that, up to machine precision,
the Riemann invariants are constant, i.e.
\begin{equation*}
	\frac{q_L^2}{2h_L^2}+g(h_L+Z_L) = \frac{q_R^2}{2h_R^2}+g(h_R+Z_R).
\end{equation*}
In addition, Neumann boundary conditions are prescribed,
and the final time is $t_{\text{end}} = 0.075$.
As discussed in \cref{sec:comments_steady},
we conjectured that such Riemann data,
although it is solution to the discrete form of Bernoulli's equation,
is not a stable steady solution.
Therefore, it should not be exactly preserved by the numerical scheme.

\begin{figure}[!ht]
	\centering
	\includegraphics{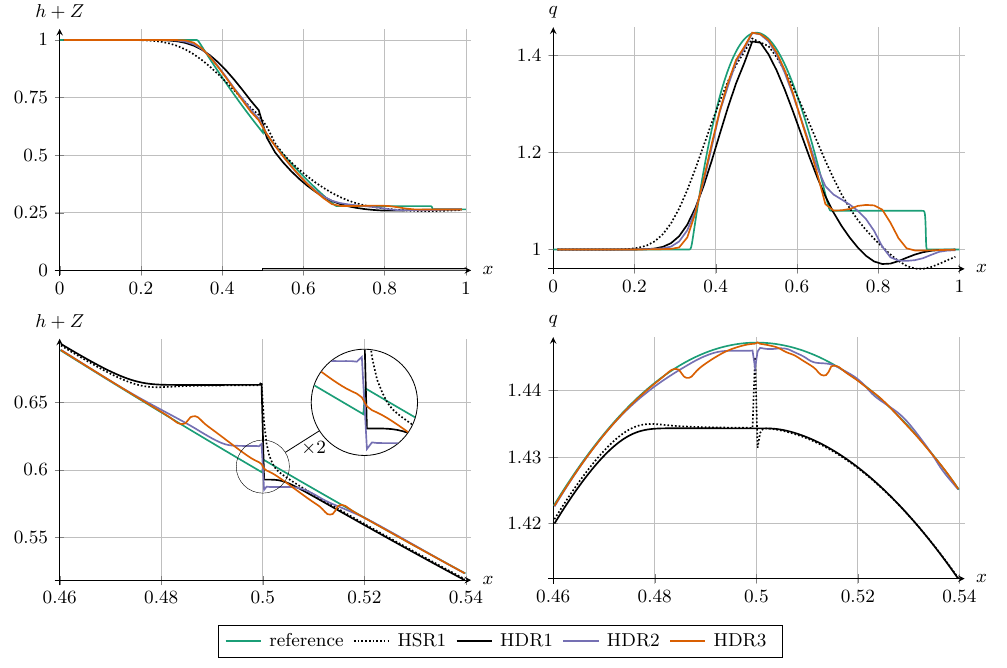}
	\caption{%
		Stationary contact from \cref{sec:stationary_contact}:
		free surface $h+Z$ (left panels)
		and discharge $q$ (right panels),
		displayed at time $t_{\text{end}} = 0.075$
		with $50$ cells (top panels)
		and $2000$ cells (bottom panels).
		The bottom panels are zoomed in on the interval (0.46, 0.54).
	}
	\label{fig:SWT_RP_small_dZ_hydro_solution}
\end{figure}

In \cref{fig:SWT_RP_small_dZ_hydro_solution},
we compare the numerical solution of the four schemes to a reference
solution computed with~$5000$ cells.
The numerical solutions with $100$ cells are displayed in the top panels,
and we observe good agreement between the numerical and reference solutions,
especially for the higher order schemes.

However, there is still a kink present around $x=0.5$,
which corresponds to the position of the topography discontinuity.
To further analyze this issue, in the bottom panels of
\cref{fig:SWT_RP_small_dZ_hydro_solution},
we provide a zoom on the interval $(0.46, 0.54)$, computed with $2000$ cells.
We observe in the bottom left panel
that the HSR$1$ and HDR$1$ schemes
display a sharp water height discontinuity around $x=0.5$.
Using the higher order HDR$2$ and HDR$3$ schemes,
the amplitude of this discontinuity decreases.
In the bottom right panel, these findings are confirmed,
although the discharge remains continuous with the HDR$1$ scheme,
contrary to the HSR$1$ scheme where a sharp oscillation is present.
Like in the case of the water height,
the higher order HDR$2$ and HDR$3$ schemes allow
a better approximation of the reference solution.

\section{Conclusion and outlook}

In this paper, we have presented an extension
\eqref{eq:hydrodynamic_reconstruction}--\eqref{eq:source}
to the hydrostatic reconstruction from \cite{AudBouBriKlePer2004}.
Applied to a numerical scheme with any consistent numerical flux function,
this hydrodynamic reconstruction possesses the following properties:
\begin{enumerate}[nosep, label=(\roman*)]
	\item consistency with the shallow water system \eqref{eq:shallow_water},
	\item preservation of moving steady solutions \eqref{eq:steady_with_constants}
	      as well as of the lake at rest,
	\item handling of transitions between wet and dry areas.
\end{enumerate}
These properties are summarized in \cref{thm:scheme_properties}.
The hydrodynamic reconstruction depends on the choice of a function $\cH$,
which has to satisfy properties \ref{item:def:H_property_dZ_0},
\ref{item:def:H_property_WB} and \ref{item:def:properties_of_H_dry_wet}.
We have exhibited such a function, and proven that is satisfies the required properties,
in \cref{prop:expression_of_H}.
Numerical experiments have confirmed that the numerical scheme
endowed with the hydrodynamic reconstruction is indeed consistent, well-balanced,
and able to treat dry/wet transitions.

Nevertheless, there are some potential improvements to the method.
First, one could design a function~$\cH$ with a more compact expression,
without losing the properties outlined in \cref{prop:expression_of_H}.
Second, one could modify the function $\cH$ to try and remove any kinks appearing
when the Froude number approaches unity.
	{\ra However, the structure of the solution is completely different
		at critical points (where $\Fr = 1$).
		Instead of the traditional Bernoulli relations,
		the slope of the water height at the critical points
		is then governed by a $p$-Laplacian-like equation.
		Changes in the nature of PDEs are widely
		recognized to pose significant numerical challenges.
		For global PDE nature changes,
		asymptotic-preserving schemes have been constructed.
		Here, the PDE nature change is local at critical points,
		and new strategies must be developed.}

\section*{Acknowledgment}
C. Berthon acknowledges the support of ANR MUFFIN ANR-19-CE46-0004.
The SHARK-FV conference has greatly contributed to this work.


\appendix

\section{\texorpdfstring{Taylor expansions of $\mathcal{H}$}{Taylor expansions of H}}
\label{sec:Taylor_expansions_H}

The goal here is to provide a Taylor expansion of the function
$\cH$ given by \eqref{eq:cH_version_2},
in the case where $\dZ > 0$, $\dh > 0$ and $1 - \Fr^2 > 0$,
when $\dZ$ approaches zero.
The computations are performed below,
where we have temporarily set $\mathbb{F} = 1 - \Fr^2$
in order to save some space.

\begin{equation*}
    \label{eq:limit_H_dZ_0}
    \begin{aligned}
        \cH
         & \underset{\dZ \to 0^+}{=}
        \frac \dh 4 \lp
        1 + \frac {\mathbb{F}} 4 \sqrt{\frac{\dh}{\dZ}}
        - \sqrt{\frac{\dh}{\dZ}} \sqrt{
            \lp \sqrt{\frac{\dZ}{\dh}} + \frac {\mathbb{F}} 4 \rp^2
            + \lp \frac{\dZ}{\dh} \rp^{\sfr 3 2}
        }
        \rp                                           \\
         & \underset{\dZ \to 0^+}{=}
        \frac \dh 4 \lp
        1 + \frac {\mathbb{F}} 4 \sqrt{\frac{\dh}{\dZ}}
        - \frac {\mathbb{F}} 4 \sqrt{\frac{\dh}{\dZ}} \sqrt{
            1
            + \frac 8 {\mathbb{F}} \lp \frac{\dZ}{\dh} \rp^{\sfr 1 2}
            + \frac {16} {\mathbb{F}^2} \frac{\dZ}{\dh}
            + \lp \frac{\dZ}{\dh} \rp^{\sfr 3 2}
        }
        \rp                                           \\
         & \underset{\dZ \to 0^+}{=}
        \frac \dh 4 \lp
        1
        + \frac {\mathbb{F}} 4 \sqrt{\frac{\dh}{\dZ}}
        \lp
        1 - \lb
        1
        + \frac 4 {\mathbb{F}} \lp \frac{\dZ}{\dh} \rp^{\sfr 1 2}
        + \frac 8 {\mathbb{F}^2} \frac{\dZ}{\dh}
        - \frac 8 {\mathbb{F}^2} \frac{\dZ}{\dh}
        + \cO \lp \dZ^{\sfr 3 2} \rp
        \rb
        \rp
        \rp                                           \\
         & \underset{\dZ \to 0^+}{=} \cO \lp \dZ \rp,
    \end{aligned}
\end{equation*}


\section{\texorpdfstring{Taylor expansions of $\mathcal{C}$}{Taylor expansions of C}}
\label{sec:Taylor_expansions_C}

In this section, we give the Taylor expansions of the function $\mathcal{C}(h_L, h_R, q_L, \dZ)$,
given by \eqref{eq:def_C_is_Fr_and_H},
when either $h_L$ or $h_R$ go to $0$.
The goal is to prove \ref{item:def:properties_of_H_dry_wet},
i.e., prove that $\mathcal{C}$ is continuous when either $h_L$ tends to $0$ and $h_R \neq 0$,
and when $h_L$ tends to $0$ and $h_R = h_L$.
Recall that
\begin{equation*}
    \mathcal{C}(h_L, h_R, q_L, \dZ) = 2 \Fr^2(h_L, h_R, q_L, \dZ) \cH(h_L, h_R, q_L, \dZ).
\end{equation*}
Since $q_L = h_L u_L$, with $u$ the velocity, we also note that
\begin{equation*}
    \label{eq:Froude_in_appendix}
    \Fr^2(h_L, h_R, h_L u_L, \dZ)
    =
    \frac{h_L^2 u_L^2 (h_L + h_R)}{2 g h_L^2 h_R^2}
    =
    \frac{u_L^2}{2 g h_R} \lp 1 + \frac{h_L}{h_R} \rp.
\end{equation*}

First, we consider the case where $h_L$ goes to $0$ and $h_R \neq 0$.
According to assumptions \eqref{eq:limh0}, in this case, $u_L$ also goes to zero.
To model this phenomenon, we assume that $u_L = u(h_L)$,
where the function $u$ is such that $u(0) = 0$.
In this case, we get, again using symbolic computation software,
\begin{equation*}
    \begin{aligned}
        \mathcal{C}
        \underset{h_L = 0^+}{=} & \pm
        \frac{\sqrt{4 g^2 h_R^2 + 16 g \sqrt{h_R \abs{\dZ}} \lp \pm u(0)^2 - 2 g (h_R \pm 2 \dZ) \rp + 4 g h_R \lp 16 g \dZ - u(0)^2 \rp + u(0)^4}}{32 g^2 \sqrt{h_R \abs{\dZ}}} \, u(0)^2 \\
        & \pm
        \frac{\sqrt{h_R} u(0)^2 - 2 g h_R \lp \sqrt{h_R} \pm 4 \sqrt{\abs{\dZ}}\rp}{32 g^2 \sqrt{h_R \abs{\dZ}}} \, u(0)^2
        +
        \cO(h_L).
    \end{aligned}
\end{equation*}
In the above expression, for the sake of clarity,
the $\pm$ symbols correspond to $\sgn(\dZ)$.
In any case, since $u(0) = 0$, this Taylor expansion shows that
\begin{equation*}
    \lim_{\substack{h_L \to 0^+ \\ h_R \neq 0}} \mathcal{C}(h_L, h_R, q_L, \dZ) = 0,
\end{equation*}
which is what we had set out to prove.

Second, we have to prove that
\begin{equation*}
    \lim_{\substack{h_L \to 0^+}} \mathcal{C}(h_L, h_L, q_L, \dZ) = 0.
\end{equation*}
However, recall from \eqref{eq:H_as_jump_times_bounded_function}
that $\cH = (h_R - h_L) \mathcal{B}$, with $\mathcal{B}$ a bounded function.
The above result is established by arguing the boundedness \eqref{eq:limh0} of the Froude number.

Therefore, property \ref{item:def:properties_of_H_dry_wet} is satisfied by $\mathcal{H}$,
when $\cH$ is given by \eqref{eq:cH_version_2}.

\bibliographystyle{plain}
\bibliography{biblio}

\end{document}